\documentclass[twoside,11pt]{amsart}

\usepackage[cmtip,all]{xy}
\usepackage{amssymb,latexsym, xspace, enumerate, amsmath, latexsym, times}
\usepackage{amssymb,amscd}

\usepackage{graphicx}
\input amssym.def
\input amssym
\input xypic
\input xy
\xyoption{all}
\def\demo{\noindent{\bf Proof. }}
\def\sqr#1#2{{\vcenter{\hrule height.#2pt
        \hbox{\vrule width.#2pt height#1pt \kern#1pt
                \vrule width.#2pt}
        \hrule height.#2pt}}}
\def\square{\mathchoice\sqr64\sqr64\sqr{4}3\sqr{3}3}
\def\QED{\hfill$\square$}

\def\m{{\mathfrak m}}
\def\n{{\mathfrak n}}

\def\p{{\mathfrak p}}

\newtheorem{Theorem}{Theorem}[section]
\newtheorem{Lemma}[Theorem]{Lemma}
\newtheorem{Corollary}[Theorem]{Corollary}
\newtheorem{Proposition}[Theorem]{Proposition}
\newtheorem{Conjecture}[Theorem]{Conjecture}

\newtheorem{Assumptions and Discussion}[Theorem]{Assumptions and Discussion}

\newtheorem{Remark}[Theorem]{Remark}
\newtheorem{Example}[Theorem]{Example}
\newtheorem{Definition}[Theorem]{Definition}
\newtheorem{Question}[Theorem]{Question}

\begin{document}
\baselineskip=16pt

\title[On the Cohen-Macaulayness of the conormal module of an ideal]
{\Large\bf On the Cohen-Macaulayness of the conormal module of an ideal}

\author[P. Mantero  and Y. Xie]
{Paolo Mantero   \and Yu Xie}

\thanks{AMS 2010 {\em Mathematics Subject Classification}.
Primary 13H05; Secondary 13H10, 13H15.}

%\thanks{The first author was partially supported by the
%NSF and the NSA}

\address{{\small Department of Mathematics, Purdue University,
West Lafayette, Indiana 47906}} \email{\small pmantero@math.purdue.edu}

\address{{\small{Department of Mathematics, University of Notre Dame,
Notre Dame, Indiana 46556}}} \email{\small yxie@nd.edu}

%\pagestyle{empty}

%\pagestyle{myheadings}

%Instead of using \demo and then \QED and \medskip we could use \begin{proof}  and \end{proof}
%\vspace{-0.1in}

\begin{abstract}
In the present paper we investigate a question stemming from  a long-standing conjecture of Vasconcelos: in a regular local ring $R$, given a perfect ideal $I$ that is a generic complete intersection, is it true that if $I/I^2$ (or $R/I^2$) is Cohen-Macaulay then $R/I$ is Gorenstein? Huneke and Ulrich, Minh and Trung, Trung and Tuan and -- very recently --  Rinaldo Terai and Yoshida, already considered this question and gave a positive answer for special classes of ideals. We give a positive answer for some classes of ideals, however, we also exhibit prime ideals in regular local rings and homogeneous level ideals in polynomial rings showing that in general the answer is negative. The homogeneous examples have been found thanks to the help of J. C. Migliore. Furthermore, the counterexamples show the sharpness of our main result. As a by-product, we exhibit several classes of Cohen-Macaulay ideals whose square is not Cohen-Macaulay. Our methods work both in the homogeneous and in the local settings.
\end{abstract}

\maketitle

\vspace{-0.2in}

For an ideal $I$ in a Noetherian local ring $R$, the first conormal module $I/I^2$  plays an important role in the process of understanding the structure of $I$. One (and the first)  main  result in this sense is a celebrated theorem of Vasconcelos stating that, if the ideal $I$ has finite projective dimension, the freeness of $I/I^2$ (as $R/I$-module) is equivalent to $I$ being a complete intersection  \cite{Va1}. Later, Vasconcelos  \cite{Va2} proposed that one could relax the assumptions in the above statement, namely he conjectured that if  the projective dimension of $I/I^2$ over $R/I$ is finite, i.e., ${\rm pd}_{R/I}(I/I^2)\,< \infty$, then $I$ must be a complete intersection. This conjecture has been proved for some classes of ideals in a Noetherian local ring by Vasconcelos (see for instance \cite{Va2} and \cite{Va3}) and, in the characteristic zero case, for any homogeneous ideal by Avramov and Herzog \cite{AH}. %but it is still wide-open in the local settings.
%On the other hand, the graded version of this conjecture has been proved by Avramov and Herzog in characteristic zero in , where they show that it holds for any ideal defining a finitely generated graded algebra over a field.
In 1987,  Vasconcelos  asked whether a similar statement could be proved in order to characterize height three Gorenstein ideals  \cite[Conjecture~B]{Va4}. This conjecture was later generalized in \cite{SVV} by removing the assumption on the height. More precisely, the conjecture asks whether the Cohen-Macaulayness of $I/I^2$, where $I$ is a {\em syzygetic} (a technical assumption) perfect ideal in a regular local ring, implies the  Gorensteiness of $I$ (\cite[Conjecture~3.12]{SVV}). The conjecture is still wide open.% (even in the original case of height $3$).% as it has been basically proved  only for almost complete intersection ideals.% and  homogeneous ideals of type two having a pure resolution \cite{Va4}.
%and suggested that a similar argument would work for any homogeneous ideal of height three having a pure resolution.

In this paper we study a question closely related to \cite[Conjecture~3.12]{SVV} (for short, we will  call it `the question'), which considers perfect ideals of any height that  are {\em generically a complete intersection} %, a class of ideals that is much wider than the one considered in the conjecture of Vasconcelos %technical condition which is, in some sense,  slightly weaker than syzygeticness
 (for detailed statement, see Question \ref{Vasc2} in Section 1).  The question can be phrased in the following slightly different way, highlighting a possible intriguing connection between Cohen-Macaulayness of the square of an ideal and the Gorenstein property. {\em Let $R$ be a regular local ring and $I$ a perfect  $R$-ideal that is generically a complete intersection.
 If $R/I^2$ is Cohen-Macaulay, does it follow that $R/I$ is Gorenstein?}.

The question is known to be true for the following classes of ideals:
\begin{itemize}
\item  perfect {\em prime} ideals of height 2 (see work of Herzog  \cite{He});
\item  licci ideals (this is one of the main results that Huneke and Ulrich proved in \cite{HU3});
\item  squarefree monomial ideals in a polynomial ring whose square is Cohen-Macaulay over any field (this is the main result in a very recent paper of Rinaldo, Terai and Yoshida \cite{RTY}).
\end{itemize}
Also, thanks to work of Minh-Trung and Trung-Tuan \cite{MT} and \cite{TT}, there is a complete description of all 2- or 3-dimensional Stanley-Reisner ideals whose second power is Cohen-Macaulay whereas the third power is not Cohen-Macaulay. In all these cases, the ideal is Gorenstein (independent of the base field), therefore providing other instances where the question is know to hold true.\\
However, in general, the question is wide open.\\
\\
In the present paper we approach the problem by means of the Hilbert function. We use it to give a positive answer to the question for some classes of ideals, namely, ideals defining stretched algebras, short algebras and algebras with low multiplicity (whose precise definitions will be given later in the paper).  Furthermore, we extend the result of Rinaldo, Terai and Yoshida and give a positive answer also for {\em any} monomial ideal (see \ref{monomial}). However, we also exhibit counterexamples and show that, in general, the question does not have a positive answer even for {\em prime} ideals.
%Our approach to the problem is by means of the Hilbert Function. It allows us to give a positive answer to the question for ideals defining stretched algebras, short algebras and algebras with low multiplicity (whose precise definitions will be given later in the paper). Furthermore, it suggested us the class of ideals where we found the counterexamples (see below).

%showing in particular Conjecture \ref{Vasc} for these classes of ideals over a regular local ring.
 %the conjecture of Vasconcelos holds true.
Combining together the main results of the first five sections, we prove the following:

\begin{Theorem}\label{Main}
Let $R$ be a regular local ring containing a field $k$ $($resp. a polynomial ring over a field $k$$)$ of characteristic $\neq 2$
 and  $I$ a  Cohen-Macaulay $($resp. homogeneous Cohen-Macaulay$)$ ideal that is generically a complete intersection. Assume that one of the following conditions holds$:$\\
$($a$)$ $I$ is a monomial ideal, under a technical assumption $($see Theorem \ref{monomial}$)$$;$\\
$($b$)$ $R/I$ is a stretched algebra$;$ \\%and ${\rm char}\,k = 0;$
$($c$)$ $R/I$ is a short algebra with socle degree at least $3;$\\
$($d$)$ $R/I$ is a short algebra with socle degree $2$ and its multiplicity satisfies some numerical conditions $($see Proposition \ref{short2}$);$\\
$($e$)$ $e(R/I)\,<\,{\rm ecodim}(R/I)+5,$ where
${\rm ecodim}(R/I)$ is the embedding codimension of $R/I.$

If $I/I^2$ $($equivalently, $R/I^2$$)$ is Cohen-Macaulay then $R/I$ is Gorenstein.
\end{Theorem}

We want to remark that parts (d) and (e) are sharp. %In particular, we exhibit a prime ideal $\p$ in a regular local ring $R$ with $e(R/\p)={\rm ecodim}(R/\p)+5$ that is a counterexample to the question.% and cannot be improved without adding extra assumptions.
In fact, the help of J. C. Migliore produced -- using CoCoa \cite{Co} --  a homogeneous level ideal $I$ defined by a set of $10$ (general) points in $\mathbb P^5$ that is a counterexample to the question. 
Remarkably, this ideal $I$ lies in the first class of ideals not covered by our main results since it defines a short algebra of socle degree $2$ (showing the sharpness of (d)) and  $e(R/I)=10={\rm ht}\,I+5$ (proving the sharpeness of (e)).

%Remarkably, this ideal $I$ lies in the first class of ideals not covered by our main results since it defines an algebra with multiplicity $e(R/I)=10={\rm ht}\,I+5$. Furthermore, the algebra defined by $I$ is short of socle degree $2$, and shows that the numerical conditions of part (d) are very sharp.

%Notice that in Theorem \ref{Main} we give a positive answer for ideals defining short algebras of socle degree at least 3. However, somewhat surprisingly, in the case of socle degree 2, we were only able to prove the same result after adding restrictions on the multiplicity. For the remaining cases, the help of J. C. Migliore produced -- using CoCoa \cite{Co} --  a homogeneous level ideal $I$ defined by a set of $10$ (general) points in $\mathbb P^5$ that is a counterexample to the question. Remarkably, this ideal $I$ lies in the first class of ideals not covered by our main results since it defines a short algebra of socle degree $2$ and  $e(R/I)=10={\rm ht}\,I+5$. %hence, Question \ref{Vasc2} has a negative answer and shows that Conjecture \ref{Vasc} cannot be extended to ideals of height at least $5$.
Using tools from linkage theory, we then deform this homogeneous ideal  to  a prime ideal $\p$
in a regular local ring $S$
with $e(S/\p)={\rm ht}\,\p+5$ such that $S/\p$ and $\p/\p^2$ are Cohen-Macaulay but $S/\p$ is not Gorenstein, proving that the answer is negative (and Theorem \ref{Main} is sharp) even for {\em prime} ideals.
This method actually shows that the whole question can be reduced to the case of prime ideals (at least in the local setting) -- see Proposition \ref{prime} and the discussion after it.

The structure of the paper is the following. In Section $1$, we prove basic facts that will be used throughout the paper, and we reduce the question  to the case of non-degenerate prime ideals. At the end of the section we explain the structure of the proof of our main results.  In Section $2$, we prove that the conormal module (equivalently,  the square) of a stretched ideal is almost never
Cohen-Macaulay and that  the question has a positive answer for stretched ideals. In Section $3$ we generalize the result of Rinaldo-Terai-Yoshida and prove that Question \ref{Vasc2} has a positive answer for monomial ideals.
In Section $4$, we prove that the conormal module (equivalently, the square) of almost every short algebra is not Cohen-Macaulay. In particular, this implies that  the question is true for almost all ideals defining short algebras.
In Section $5$, we employ results from Sections 2 and 4 to show that  the question   has a positive answer for ideals defining algebras whose multiplicity is at most the embedding codimension plus 4. We then provide an application to algebroid curves.
In Section $6$, we exhibit a set of $10$ (general) points in $\mathbb P^5$ (found by J. C. Migliore)
showing that, in general, the question has a negative answer in the homogeneous setting. We also prove that the local version of the question has a negative answer even for prime ideals and show the sharpness of Theorem \ref{Main}. %and of results presented in Section $3$.
%We also provide some counterexamples in higher dimensions and
%The embedding codimension of these examples is 5.
We then conjecture that,  for {\em any} embedding codimension $c\,\geq\,5$, the homogeneous ideal $I$ defined by a set of $c+1+\left\lceil\frac{c(c-1)}{6}\right\rceil$ general points in $\mathbb P^c$ has the property that $I/I^2$ is Cohen-Macaulay, but $R/I$ is not Gorenstein.
\bigskip

\section{Preliminaries}

In this section we fix some notation and prove basic facts that will be used throughout the paper. In particular, we show that in Question \ref{Vasc2}  one can always assume the ideal $I$ to be `non-degenerate', i.e., $I\subseteq \m^2$ (see Proposition \ref{NonDeg2}), and prime (see Proposition \ref{prime}).

Recall that in a Noetherian  ring $R$, an ideal $I$ is {\em generically a complete intersection} if $I_{\p}$ is a complete intersection $R_{\p}$-ideal for every $\p \in {\rm Ass}_{R}(R/I)$, the set of  associated prime ideals of $R/I$.  We are now able to state Vasconcelos's question in the form that it appears in \cite{HU3}.

%We now state the generalization of Conjecture \ref{Vasc}
\begin{Question}$($\cite{HU3}$)$\label{Vasc2}
Let $R$ be a regular local ring $($resp. a polynomial ring over a field$)$ and $I$ a Cohen-Macaulay  $($resp. homogeneous Cohen-Macaulay$)$ ideal that is generically a complete intersection. If $I/I^2$ is Cohen-Macaulay then does $R/I$ have to be Gorenstein?
\end{Question}

Notice that from the Cohen-Macaulayness of $I$ and the short exact sequence
$$0 \longrightarrow I/I^2 \longrightarrow R/I^2 \longrightarrow R/I \longrightarrow 0$$
it follows that $I/I^2$ is Cohen-Macaulay if and only if $R/I^2$ is. Hence one can phrase Question \ref{Vasc2} in the following slightly different way, involving only $I$ and $I^2$.

\begin{Question}\label{Vasc3}
Let $R$ be a regular local ring $($resp. a polynomial ring over a field$)$ and $I$ a generically complete intersection ideal. If both $R/I$ and $R/I^2$ are Cohen-Macaulay  $($resp. homogeneous Cohen-Macaulay$)$  then does it follow that $R/I$ is Gorenstein?
\end{Question}

The classes of ideals for which Question \ref{Vasc3} has a positive answer, reflect an interesting connection between the Cohen-Macaulayness of $I^2$ and the Gorenstein property of $I$.

%we gather some results that let us consider some reductions, in particular we show that without loss of generality one may assume (both for the conjecture of Vasconcelos and the question generalizing it) that the ideal and  prime, see \ref{prime} and \ref{NonDeg2}. We also prove Proposition \ref{Generalized} that will be used often in the rest of the paper.

%We begin by stating Vasconcelos's Conjecture,
%\begin{Conjecture}$($\cite[Conjecture~B]{Va4}$)$\label{Vasc}
%Let $(R,\m)$ be a regular local ring and $I$ a height three Cohen-Macaulay ideal that is syzygetic. If $I/I^2$ is Cohen-Macaulay then $R/I$ has to be Gorenstein.
%\end{Conjecture}

%We will prove that for stretched algebras (see \ref{defstretched}), almost every short algebra (see \ref{defshort}, \ref{short} and \ref{short2}) and algebras having low multiplicity (see \ref{LowMult}),
% the answer is positive, showing in particular that Vasconcelos Conjecture holds true for these classes of ideals. However, in general the answer to Question \ref{Vasc2} is negative, even if one considers the analogous question in the graded settings.

Since Questions \ref{Vasc2} and \ref{Vasc3} are equivalent, from now on we will only refer to Question \ref{Vasc2}.
To study this question, we first reduce it to the `non-degenerate' case, i.e.,  $I\subseteq \m^2$,  where $\m$ is the maximal ideal of $R$. This is accomplished in Proposition \ref{NonDeg2}. In order to prove it, we first need to collect a couple of lemmas. The first one gives a condition for the conormal module $I/I^2$ to contain a free summand as $R/I$-module.
%The remaining part of this section is devoted both to reduce the problem to the case where $I\subseteq \m^2$ (the so called {\em non-degenerate case}) and to relate the multiplicity of $R/I^2$ with $e(R/I)$.\\
%Next two results serve for reducing to the non-degenerate case.

\begin{Lemma}\label{splitting}
Let $R$ be a Noetherian  ring and  $I=(\underline{x})+J$, where  $J$  is an $R$-ideal and $\underline x=x_1,\dots,x_t$ \,forms a regular sequence on $R/J$. Then $I/I^2$ has a free summand $($as $R/I$-module$)$ of rank $t$. In fact,  $(\underline x+I^2)/I^2\simeq (R/I)^t$ and $I/I^2= (\underline x+I^2)/I^2 \oplus K$ for some submodule $K$ of $I/I^2$.
\end{Lemma}

\demo
If $I=R$ then we are done. So we may assume that $I$ is a proper $R$-ideal.
Set $\overline{R}=R/J$ and denote by $^-$ the images in $\overline{R}$, then ${\overline I}=(\overline{x_1},\dots,\overline{x_t})$ is a complete intersection $\overline R$-ideal. Hence $\overline{I}/\overline{I}^2$ is a free $\overline{R}/\overline{I}$-module of rank $t$.
Now the first part of the lemma follows by the exact sequence $I/I^2\rightarrow I/(J+I^2) \rightarrow 0$ and  the fact that $I/(J+ I^2)\simeq (I/J)/(I/J)^2\simeq \overline{I}/\overline{I}^2\simeq (\overline{R}/\overline{I})^t\simeq(R/I)^t$.
Finally since $(\underline x+I^2)/I^2\simeq I/(J+\!I^2)\simeq (R/I)^t$,
the second part of the lemma follows from the split exact sequence  $0 \rightarrow K \rightarrow I/I^2\rightarrow (\underline x+I^2)/I^2\rightarrow 0$.
%and let $K$ be the kernel of this map.  Since $(\underline x+I^2)\simeq (R/I)^t$, this proves that $I/I^2= (\underline x+I^2) \oplus K$, thus concluding the proof.
\QED

\bigskip

%Recall in a  Noetherian local ring $(R,\m)$,  $R^*={\rm gr}_{\m}(R)=\bigoplus_{i\,\geq\,0}\m^i/\m^{i+1}$ is the {\em associated graded ring} of $R$ with respect to $\m$.
% For an element $0\neq a\in R$, there is a unique $i$ with $a \in \m^i\backslash\m^{i+1}$. Then the {\em initial form} of $a$ is: $a^*=(a+\m^{i+1})/\m^{i+1} \in R^*$ and the {\em initial degree} ${\rm indeg}(a)$ of $a$  is $i$, i.e., the degree of $a^*$ in the graded ring $R^*$. %Similarly,   the {\em ideal of initial forms} of  an $R$-ideal $J$ is defined as  $J^*=( a^* \,|\, a\in J)R^*$. It is well known that ${\rm gr}_{\m}(R/J)\simeq R^*/J^*$.
%The following lemma is a well-known result showing that if the initial form of an element is regular modulo the ideal $J^*$,  then the element itself is regular modulo $J$.

%\begin{Lemma}\label{initial}
%Let $R$ be a Noetherian local ring and $J$ an $R$-ideal.  Let $x$ be an element of $R$ such that $x^*$ is regular both on $R^*$ and on $R^*/J^*$. Then $x$ is regular on $R/J$.
%\end{Lemma}

%We now employ the previous two lemmas to obtain the following Proposition \ref{NonDeg}.
Recall that a Noetherian local ring $R$ is said to be {\em equicharacteristic} if it contains a field. $R$ has {\em mixed characteristic} if it does not contain a field. We are now able to prove the following result.

\begin{Proposition}\label{NonDeg}
Let $(R,\m)$ be a regular local ring containing a field and  $I$ an $R$-ideal containing an element $x\in \m\backslash \m^2$. Then $I/I^2$ contains a free summand $($as $R/I$-module$)$ of rank $1$.
\end{Proposition}

\demo The statement is trivial if $I=R$, hence we may assume that $I$ is a proper $R$-ideal. We first prove that one can always reduce to the case where $R$ is complete.  Set $S=R/I$, $M=I/I^2$, let $N=(xR+I^2)/ I^2 =\bar{x}S\subseteq M$ be the cyclic $S$-submodule of $I/I^2$ generated by the image of $x$ and assume that $\widehat{N}$ is a free direct summand of $\widehat{M}$ of rank $1$. We need to show that $N$ is a free direct summand of $M$ of rank $1$.
Let $K$ be the kernel of the presentation map $\phi :S \rightarrow N \rightarrow 0$ defined by $\phi(1_S)=\bar{x}$. Since the induced map $\hat{\phi}:\widehat{S}\rightarrow \widehat{N}\rightarrow 0$ is an isomorphism of $\widehat{S}$-modules, it follows that $K=0$, proving the freeness of $N$. Next, notice that, by assumption,  we have the epimorphism of $\widehat{S}$-modules $\hat{j}: \widehat{M}\rightarrow \widehat{N}$ with $\hat{j}(\bar{x})=\bar{x}$ and $\hat{j}(\tilde{f}_i)=0$, $1\leq i \leq r$, where $\bar{x},  \tilde{f}_1, \ldots, \tilde{f}_r$ form a set of minimal generators of $\widehat{M}$.
Since ${\rm Hom}_{\widehat{S}}(\widehat{M}, \widehat{N})\simeq \widehat{{\rm Hom}_S(M,N)}$, there exists a homomorphism of $S$-modules $\pi: M \rightarrow N$ such that $\hat{j}-\hat{\pi}\in \m \widehat{{\rm Hom}_S(M,N)}\simeq \m {\rm Hom}_{\widehat{S}}(\widehat{M}, \widehat{N})$, where $\hat{\pi}:\widehat{M}\rightarrow \widehat{N}$ denotes the map induced by $\pi$.
Hence $(\hat{j}-\hat{\pi}) (\widehat{M})\subseteq \m \widehat{N}$ and thus $\hat{\pi}(\widehat{M})+\m \widehat{N}=\hat{j}(\widehat{M})=\widehat{N}.$ By Nakayama's Lemma, $\hat{\pi}(\widehat{M})=\widehat{N}$, i.e., $\hat{\pi}: \widehat{M}\rightarrow \widehat{N}$ is surjective. Since $\widehat{S}$ is faithfully flat over $S$, the map $\pi: M\rightarrow N$ is also surjective. Therefore $\pi(\bar{x})=\mu \bar{x}$ for some unit $\mu$ in $S$. If one denotes by $\iota$ the injective map $ N \hookrightarrow M$ defined by $\iota(\bar{x})=\bar{x}/\mu$,  then we have $\pi\circ \iota \equiv {\rm id}_N$,  and hence $N$ is a free direct summand of $M$.
%Since ${\rm Hom}_{\hat{S}}(\hat{M}, \hat{N})\simeq \widehat{{\rm Hom}_S(M,N)}$, there exists a homomorphism of $S$-modules $j: M \rightarrow N$ such that $\hat{j}-j\in \m \widehat{{\rm Hom}_S(M,N)}\simeq \m {\rm Hom}_{\hat{S}}(\hat{M}, \hat{N})$.
%Hence $(\hat{j}-j) (\hat{M})\subseteq \m \hat{N}$ and thus $j(\hat{M})+\m \hat{N}=\hat{j}(\hat{M})=\hat{N}.$ By Nakayama's Lemma, $j(\hat{M})=\hat{N}$, i.e., $j: \hat{M}\rightarrow \hat{N}$ is surjective. Since $\hat{S}$ is faithfully flat over $S$, $j: M\rightarrow N$ is also surjective. If one denotes by $i$ the embedding map $i: N \hookrightarrow M$, since $j(x)=\hat{j}(x)=x$ then we have $j\circ i \equiv {\rm id}_N$,  and $N$ is a free direct summand of $M$.

We may then assume that $R$ is a complete regular local ring containing a field $k$. Write $R\simeq k[\![x,x_1,\dots,x_s]\!]$ for some $s$ and $I=(x,J)$ for some $R$-ideal $J=(f_1,\dots,f_r)$. %Since each $f_i$ is a power series in the variables $x,x_1,\dots,x_s$ with coefficients in $k$,
Now use $x$ to write $I=(x,J')$ with $J'=(f_1',\dots,f_r')$ and each $f_i'$ is  in $k[\![x_1,\dots,x_s]\!]$. Then $x$ is clearly regular on $R/J'$ and by Lemma \ref{splitting}, $I/I^2$ contains a free direct summand of rank $1$, finishing the proof.
\QED

\bigskip
Notice that the assumption of $R$ containing a field cannot be removed. In \cite{Ho}, Hochster constructed an example to provide an obstruction to Grothendieck's Lifting Problem. We employ it here to show that the assumption of $R$ being equicharacteristic is needed: if we assume $R$ to have mixed characteristic, Proposition \ref{NonDeg} fails to be true even in the unramified case. %In fact, in the mixed characteristic case the above Proposition does not hold.

\begin{Example}$($\cite[Example~1]{Ho}$)$
There exists an unramified complete regular local ring $(A,\n)$ and an $A$-ideal $I$ with $x\in I$ and $x \notin \n^2$ such that the $A/I$-submodule $N$ of $I/I^2$ generated by the image of $x$ has non zero annihilator. %(as $A/I$-module).
In particular, $N$ is not even a free $A/I$-submodule of $I/I^2$.
\end{Example}

\demo
Let $V$ be a complete discrete valuation ring with maximal ideal generated by $2$, and let $A=V[\![X_1,X_2,Y_1,Y_2,Z_1,Z_2]\!]$ be a power series ring in $6$ analytic variables over $V$. Set $q=X_1X_2+Y_1Y_2+Z_1Z_2$, $x=2$ and $I=(2,q,X_1^2,X_2^2,Y_1^2,Y_2^2,Z_1^2,Z_2^2)A$. We show that $I^2:2A \nsubseteq I$. In fact, set $D=X_1X_2Y_1Y_2+Y_1Y_2Z_1Z_2 + X_1X_2Z_1Z_2$, then $2D = q^2 - X_1^2X_2^2 - Y_1^2Y_2^2 - Z_1^2Z_2^2 \in I^2$ (since all the summands are in $I^2$). To show that $D \notin I$, set $K=V/2V$, $H= (2,X_1^2,X_2^2,Y_1^2,Y_2^2,Z_1^2,Z_2^2)A$ and denote by $'$ residues modulo $H$. Then $I'$ is a principal ideal generated by $q'$ and any degree 4 element in $I'$ is in the $K$-span of $q'f'$ where $f \in V$ runs over all square-free monomials of degree $4$. It is not hard to see that $D'$ is not in this span, proving that $D\notin I$.
\QED

\bigskip

Similar examples can be produced in any positive characteristic, showing that Proposition \ref{NonDeg} does not hold true in the mixed characteristic case.

As a result of Proposition \ref{NonDeg} and Lemma \ref{splitting},  we can now reduce Vasconcelos's question  to the non-degenerate case.

\begin{Proposition}[Reduction to the non-degenerate case]\label{NonDeg2}
Let $(R,\m)$ be a regular local ring containing a field and  $I$ an $R$-ideal containing an element $x \in \m \backslash \m^2$. Set $\overline{R}=R/(x)$ and $\overline{I}=I\overline{R}$. Then Question \ref{Vasc2} holds true for $I$  if and only if it holds true for $\overline{I}$.
\end{Proposition}

\demo
Since $R/I\simeq \overline{R}/\overline{I}$,  $R/I$ is Cohen-Macaulay or Gorenstein if and only if $\overline{R}/\overline{I}$ is.
It is also easy to see that  $I$ is generically a complete intersection if and only if $\overline{I}$ is.
Also, by Proposition \ref{NonDeg} and Lemma \ref{splitting}, $(x+I^2)/I^2\simeq R/I$ and $I/I^2= (x+I^2) /I^2\oplus K$ for some submodule $K$ of $I/I^2$.
Therefore $\overline{I}/\overline{I}^2\simeq K$ and it follows that $I/I^2$ is Cohen-Macaulay if and only if $\overline{I}/\overline{I}^2$ is.
% where $K$ is the kernel of the map $\Phi$ defined in Lemma \ref{splitting}.
%Then there exists a short exact sequence
\QED
\bigskip

%NEW VERSION
%\begin{Proposition}[Reduction to the non-degenerate case]\label{NonDeg2}
%Let $(R,\m)$ be a regular local ring  and  $I$ an $R$-ideal containing an element $x \in \m \backslash \m^2$. Set $\overline{R}=R/(x)$ and $\overline{I}=I\overline{R}$.  Then\\
%(a) $I$ is Cohen-Macaulay $($Gorenstein, respectively$)$ if and only if $\overline{I}$ is.\\
%(b) $I/I^2$ is Cohen-Macaulay if and only if $\overline{I}/\overline{I}^2$ is.

%Question \ref{Vasc2} holds true for $I$  if and only if it holds true for $\overline{I}$.
%\end{Proposition}

%\demo
%(a) Since $R/I\simeq \overline{R}/\overline{I}$,  $R/I$ is Cohen-Macaulay or Gorenstein if and only if $\overline{R}/\overline{I}$ is.\\
%(b) By Proposition \ref{NonDeg} and Lemma \ref{splitting}, $(x+I^2)/I^2\simeq R/I$ and $I/I^2= (x+I^2) /I^2\oplus K$ for some submodule $K$ of $I/I^2$.
%Therefore $\overline{I}/\overline{I}^2\simeq K$ and it follows that $I/I^2$ is Cohen-Macaulay if and only if $\overline{I}/\overline{I}^2$ is.
% where $K$ is the kernel of the map $\Phi$ defined in Lemma \ref{splitting}.
%Then there exists a short exact sequence
%\QED

%\bigskip

Recall that for an ideal $I$ in a regular  local ring $R$, the {\rm embedding codimension} of $\overline{R}=R/I$  is defined to be ${\rm ecodim}(\overline{R})=\mu(\overline{\m})-{\rm dim}\,\overline{R}$, where $\mu(\overline{\m})$ denotes the minimal number of generators of the maximal ideal $\overline{\m}$ of $\overline{R}$. Notice that in general the height of $I$ is greater than or equal to the embedding codimension, and if   $I$ is non-degenerate, i.e., $I\subseteq \m^2$, then ${\rm ht}\,I={\rm ecodim}(R/I)$.

The following lemma, consequence of work of Huneke and Ulrich, gives a positive answer to Question \ref{Vasc2} for all ideals with ${\rm ecodim}(R/I)\leq 2$. Therefore, throughout the paper we will always deal with the case where the embedding codimension is at least $3$.
%The following remark shows that the condition $c\geq 3$ is harmless, since Question \ref{Vasc2} is known for $c\leq 2$ as one can use the linkage to see it.
%MAYBE HERE WE NEED TO ADD THE ASSUMPTION OF R CONTAINING A FIELD? OR CAN WE FIX IT?
\begin{Lemma}\label{LC}$($\cite[2.8]{HU3}$)$
%Let $R$ be either a regular local ring or a polynomial ring over a field. Let $I$ be as in Question \ref{Vasc2}. If then
Question \ref{Vasc2} has a positive answer if ${\rm ecodim}(R/I)\leq 2$ and either $R$ contains a field or $I\subseteq \m^2$.
\end{Lemma}

\demo
Let $\m$ be the maximal ideal or the homogeneous maximal ideal of $R$.
By Proposition \ref{NonDeg2}, one can assume $I\subseteq \m^2$. Then ${\rm ht}\,I={\rm ecodim}(R/I) \leq 2$.
If ${\rm ht}\, I=1$, it is easy to see that $I$ is a licci ideal.
If ${\rm ht}\, I=2$, $I$ is licci by a classical result of Ap$\grave{{\rm e}}$ry-Ga\`eta (see \cite{Ap}, \cite{Ga}). In both cases, \cite[2.8]{HU3}  implies that Question \ref{Vasc2} is true for $I$.
\QED

\bigskip

 %The following remark shows that the condition $c\geq 3$ is harmless, since Question \ref{Vasc2} is known for $c\leq 2$ as one can use the linkage to see it.
%Hence, we will always deal with the case where  the embedding codimension is at least $3$.

Next, we want to show that Question \ref{Vasc2} can be reduced to the case of prime ideals.
%We should also mention that the last part of Theorem \ref{Main} follows from a more general statement: without loss of generality one can use linkage to reduce Conjecture \ref{Vasc} and Question \ref{Vasc2} to the case where $I$ is a prime ideal.
To do it, we first recall some definitions from linkage theory. For more details we recommend \cite{HU2}.

Let $(S,J)$ and $(R,I)$ be pairs, where  $(S,\n)$ and $(R,\m)$ are Noetherian local rings and
$J\subseteq S$ and $I\subseteq R$ are ideals. %one has the following definitions: \\
One says that $(S,J)$ and  $(R,I)$ are {\em isomorphic}, written $(S,J)\simeq (R,I)$, if there is an isomorphism of rings  $\phi: S \rightarrow R$ such that $\phi(J)=I$.
$(S,J)$ is called a {\em deformation} of $(R,I)$ if there is a sequence $\underline x=x_1,\dots,x_t$ in $S$ regular both on  $S$ and on $S/J$ such that $(S/(\underline x), J+(\underline x)/(\underline x))\simeq (R,I)$.
 Let $\underline Z=Z_1,\dots,Z_t$ be variables over $R$,\,  $\underline a=a_1,\dots,a_t$ elements in $R$ and $\underline Z-\underline a =Z_1-a_1,\dots,Z_t-a_t$. Let $J\subseteq R[\underline Z]_{(\m,\, \underline Z-\underline a)}$ be an ideal. One says that $(R[\underline Z]_{(\m, \,\underline Z-\underline a)}, J)$ is a {\em specialization} of $(R,I)$ if $(R[\underline Z]_{(\m,\,\underline Z-\underline a)}, J)$ is a deformation of $(R,I)$ via the regular sequence $\underline Z-\underline a$.

Let $R$ be a Gorenstein ring and $I=(f_1,\dots , f_n)$  an unmixed $R$-ideal of grade $g>0$. Let $X$ be a generic
$g \times n$ matrix of variables  over $R$ and  $\underline{\alpha}=\alpha_1,\dots\alpha_g$  the  regular sequence in $R[X]$ defined by $(\underline{\alpha})^t = X (\underline f)^t$, where $\underline f=f_1,\dots , f_n$.
Then $L_1(\underline f)$=$(\underline{\alpha})R[X]:_{R[X]} IR[X]$ is the {\em first generic link} of $I$.
Since $L_1(\underline f)$ is independent
of the generating set  $\underline f=f_1,\dots , f_n,$ up to isomorphism of pairs (after adjoining variables) \cite[2.11]{HU2},
we then write $L_1(I) = L_1(\underline f)$ and by iteration of this process
one defines the {\em $i$-th generic link} $L_i(I)$ as $L_i(I)=L_1(L_{i-1}(I))$ for any $i \geq 1$.

% proved much more. In fact, the same proof gives the following lemma showing that without loss of generality one can add the important assumption that $I$ is a prime ideal. More precisely:
%we are able to show the following important remark that one can always assume the ideal to be prime. This statement will follow easily from the next Lemma, where we invoke Linkage Theory to deform any counterexample to Vasconcelos Conjecture to make it become a prime ideal $\p$ in a regular local ring that is still a counterexample to Conjecture \ref{Vasc}.
Using tools from linkage, we now prove that, when dealing with Question \ref{Vasc2}, one can always work with {\em prime} ideals instead of generically  complete intersection ideals.
More precisely:

\begin{Proposition}\label{prime}
Any generically  complete intersection ideal that gives a negative answer to Question \ref{Vasc2}\, admits a deformation that is a prime ideal of the same height as the original ideal and gives a negative answer to Question \ref{Vasc2} too.
\end{Proposition}

\demo
Let $I$ be generically a complete intersection with  $R/I$  and $I/I^2$   Cohen-Macaulay, and assume $R/I$ is not Gorenstein. %We will show that there exists a prime ideal $\p$ in a (different) regular local ring $S$ with $S/\p$ and $\p/\p^2$ Cohen-Macaulay, but $S/\p$ is not Gorenstein.
Let $L_2(I)\subseteq R[\underline{Z}]$ be the second generic link of $I$ in a polynomial extension of $R$.  By  \cite{HU2},  there exists a sequence $\underline{a}$ of $R$ such that if one sets $S=R[\underline{Z}]_{(\m, \underline{Z}-\underline{a})}$ and $\p=L_2(I)S$, then $\p$ is a prime ideal of the same height as $I$ and $(S, \p)$ is a specialization  of $(R,I)$. By \cite[2.1]{Hu}, one has that $\p/\p^2$ is also Cohen-Macaulay. Since $(S,\p)$ is a specialization of $(R,I)$, it is easy to see  that $S/\p$ is Cohen-Macaulay but not Gorenstein.
\QED

\bigskip

As a consequence, if one could prove that Question \ref{Vasc2} has a positive answer for any prime ideal then, by Proposition \ref{prime}, the question would automatically be settled in full generality. At the same time, if there exists a counterexample to the question, by \ref{prime} there exists also a prime ideal that is a counterexample. Therefore, one can always assume that the ideal $I$ in Question \ref{Vasc2} is {\em prime}.

 %implies that one needs only to focus on prime ideals to give any answer to Question \ref{Vasc2}.
However, in the following we will not use the extra assumption that $I$ is a prime ideal. In fact, in this paper we prove results regarding the conormal module (or the square) of some classes of ideals and it seems natural to present them in the most general context, that is for generically complete intersection ideals.
Proposition \ref{prime} will be employed in the last section to produce a prime ideal providing a negative answer to the general version of Question \ref{Vasc2}.

%From the proof of Proposition \ref{NonDeg} it follows that $I=(y,J)$ for some $R$-ideal $J$ with $y$ regular on $R/J$.

%\begin{Remark}\label{homog}
%Let $(R,\m)$ be a regular *local ring, $x_1,\dots,x_c$ a homogeneous regular system of parameters of $\m$, $I$ a homogeneous $\m$-primary ideal. If $I$ is generated by homogeneous elements $f_1,\dots,f_r$ where each $f_i$ is a combination of monomials in $x_1,\dots,x_c$ of the same degree $d_i$ with coefficients that are units in $R$, then $R/I$ is canonically graded (in the sense of \cite{E}), therefore $(I^*)^n=(I^n)^*$ for any $n\,\geq\,1$. In particular $\lambda(R/I^2)\,>\,(c+1)\lambda(R/I)$ if and only if $\lambda(R^*/(I^*)^2)\,>\,(c+1)\lambda(R^*/I^*)$.\\
%\end{Remark}

\medskip
The last result of  this section is a useful formula for computing the multiplicity of the
square of any generically  complete intersection ideal in terms of the multiplicity of the ideal itself. %The formula is slightly different from the one usually found in literature because it uses the embedding codimension.

\begin{Proposition}\label{Generalized}
Let $R$ be an equidimensional catenary  Noetherian  local ring and $I$ an ideal that is generically a complete intersection of height $c$. Then $e(R/I^{2})=(c+1)e(R/I)$.
\end{Proposition}

\demo %Let $\m$ be the maximal ideal or the homogeneous maximal ideal of $R$.
By the Associativity formula, one has
$$e(R/I^{2})=\sum_{P}\lambda(R_P/I_P^{2})e(R/P),$$
%\lambda(R/(I^{(2)}, \underline x))=
where $P$ ranges over all minimal primes of $I$ of maximal dimension.
%Observe that for every $P \in {\rm Min}(I)$,  $I_P^{(2)}=I_P^2$.
Since $R$ is equidimensional and catenary,\,  $I_P$ is a complete intersection of height $c$ and hence $I_P/I_P^2$ is a free $R_P/I_P$ module of rank $c$. Now consider the short exact sequence
$$0 \longrightarrow I_P/I_P^2 \longrightarrow R_P/I_P^2 \longrightarrow R_p/I_P \longrightarrow 0,$$
by additivity of lengths, we obtain that $\lambda(R_P/I_P^2)=(c+1)\lambda(R_P/I_P)$.
Therefore
$$\begin{array}{ccl}
e(R/I^{2})&=& \sum_{P}\lambda(R_P/I_P^2)e(R/P)\\
&=&\sum_{P}(c+1)\lambda(R_P/I_P)e(R/P)\\
&=&(c+1)\sum_P\lambda(R_P/I_P)e(R/P)\\
&=&(c+1)e(R/I)
\end{array}$$
where the last equality follows again by the Associativity formula.
\QED

\bigskip

Proposition \ref{Generalized} will be used frequently to prove the
 main results of the next several sections. Also, since some of the proofs of these results may look technical, we want to sketch here the underlying ideas. The main steps are:
\begin{enumerate}[I.]
\item using the Macaulayness of $R/I$ and $R/I^2$, and that $I$ is generically a complete intersection, reduce Question \ref{Vasc2} to proving that $e(R/I^2)=(c+1)e(R/I)$, where $c={\rm ecodim}(R/I)$.
\item After going modulo a (special) minimal reduction of the maximal ideal of $R/I$, reduce to the case where $I$ is $\m$-primary ($I$ may not be generically a complete intersection anymore, but we don't need it in the rest of the proof).
\item Use the Hilbert function of $R/I$ to estimate the Hilbert function of $R/I^2$ -- this step is the longest and hardest, especially in the local setting. It requires the assumptions on the structure of $R/I$ (e.g. stretched, short, etc.) and involves several {\em ad hoc} methods and strategies.
\item Use step III. to estimate the multiplicity of $R/I^2$ in terms of the multiplicity of $R/I$ and show that it leads to a contradiction with the equality found in step I.
\end{enumerate}

\bigskip

%Also, since some of the proofs of these results may look technical, we would like to sketch the underlying ideas. If
% $I/I^2$ is Cohen-Macaulay then $R/I^2$ is also Cohen-Macaulay and there exists an $R$-ideal $J$ whose images in $R/I$ and $R/I^2$, respectively, give a minimal reduction of the maximal ideals of $R/I$ and $R/I^2$. If one denotes by $^{-\!\!-}$ the images in $\overline{R}=R/J$, then  $e(R/I)=e(\overline{R}/\overline{I})$ and $e(R/I^2)=e(\overline{R}/\overline{I}^2)$.
%Our goal is to estimate the Hilbert function of $\overline{R}/\overline{I}^2$ by means of the Hilbert function of $\overline{R}/\overline{I}$ in order to estimate $e(\overline{R}/\overline{I}^2)$ and then get a contradiction.
%Usually, in the local setting, giving tight estimates for $\overline{R}/\overline{I}^2$ based on the Hilbert function of $\overline{R}/\overline{I}$ is a very hard task. In this work, we accomplish this making computations on the Hilbert function of $\overline{R}/\overline{I}^2$  that use the structure of $\overline{R}/\overline{I}$ given by our assumptions on $R/I$ (e.g. $R/I$ stretched, short, etc.).
%We then employ these estimates, sometimes together with the non-Gorenstein assumption on $R/I$, to prove that the multiplicity of $R/I^2$ is different from $(c+1)e(R/I)$, contradicting the equality of Proposition \ref{Generalized}. This shows that if $R/I$ is not Gorenstein then $I/I^2$
%(or, equivalently, $I^2$) cannot be Cohen-Macaulay.

We would also like to remark that in this paper we mainly focus on the local case, since it is the most problematic one. In fact, all our proofs work also in the homogeneous settings, and for some results (e.g. stretched algebra case or the reduction to the non-degenerate case) the homogeneous setting actually allows one to  give much simpler proofs.

\bigskip

\section{Stretched algebras}

In this section, we study the conormal module of ideals defining stretched algebras and show that it is almost never Cohen-Macaulay. We start by  setting up some notation.

Let $(A,\n)$ be an Artinian local ring. We recall  that the {\em socle degree} of $A$, ${\rm socdeg}(A)$, is the (only) positive integer $s$ with $\n^{s+1}=0$ and $\n^s\neq 0$. Notice that with this definition, ${\rm socdeg}(A)$ is exactly one unit smaller than the Loewy length of $A$. The {\it socle} of $A$, ${\rm soc}(A)=0:_A \n$, is a subset of $A$ consisting of the elements of $ A$ annihilated by the maximal ideal $\n$. 
${\rm Soc}(A)$ is  a non-trivial $A/\n$-vector space whose dimension, $\tau(A)={\rm dim}_{A/\n}\,{\rm soc}(A)$, is called the ({\em Cohen-Macaulay}) {\em type} of $A$. In general if  $(R,\m)$ is a  Cohen-Macaualy local ring with $|R/\m|=\infty$, the {\em socle degree}, ${\rm socdeg}(R)$, the {\it socle}, ${\rm soc}(R)$, and the {\em type},  $\tau(R)$, are defined respectively to be  ${\rm socdeg}(R/J)$, ${\rm soc}(R/J)$ and $\tau(R/J)$, where $J$  is a minimal reduction of $\m$ generated by a system of parameters. The ring $R$ is {\em Gorenstein} if in addition $\tau(R)=1$. Finally for any Noetherian local  ring $R$, we will use $e=e(R)$ to denote the Hilbert-Samuel multiplicity of  $R$ with respect to its maximal ideal $\m$.

%\begin{Notation}
%Let $(R,\m)$ be  a Noetherian local  ring and $I$ an $R$-ideal. We always use $e=e(R/I)$ for the Hilbert-Samuel multiplicity of $R/I$ with respect to the maximal ideal $\m$, $c={\rm ecodim}(R/I)$ for the embedding codimension of $R/I$ and $s=s(A)={\rm socdeg}(R/I)$ for the socle degree of $R/I$.
%Also, we use $R^*={\rm gr}_{\m}(R)$ for the associated graded ring with respect to the maximal ideal $\m$, $I^*$ for the initial ideal of $I$ and  $a^*$ for the initial form of $a\in R.$
%\end{Notation}

We now recall the definition of stretched algebras which is one of the main classes of ideals that we will study in this paper.

\begin{Definition}\label{defstretched}
An Artinian local ring $(A,\n)$ is said to be {\rm stretched} if $\n^2$ is a principal ideal.
In general, for a Cohen-Macaulay local ring $(R,\m)$  with $|R/\m|=\infty$, we say  that $R$ is {\rm stretched}  if
$R/J$ is stretched for any $J$ minimal reduction of $\m$ generated by a system of parameters.
A Cohen-Macaulay ideal $I$ in a Noetherian local ring $R$ is {\rm stretched} if $R/I$ is.
\end{Definition}

An immediate consequence of the definition is that one can easily describe the Hilbert function $HF_A$ of a stretched  Artinian local ring $(A,\n)$:
$$1 \qquad c \qquad 1 \qquad 1 \qquad \ldots \qquad 1 \qquad 0_{\longrightarrow}$$
In fact, if $A$ has the above Hilbert function then $\mu(\n^2)=1$, which implies that $A$ is stretched. On the other hand, if $\mu(\n^2)=1$, then $HF_A(2)=HF_{{\rm gr}_{\n}(A)}(2)=1$. Now Macaulay's bound applied to ${\rm gr}_{\n}(A)$ shows that $HF_A(i)=HF_{{\rm gr}_{\n}(A)}(i)$ is at most $1$ for any $i\geq 2$, proving that $A$ has the above Hilbert function.

Another easy consequence is that in a stretched Artinian  algebra $(A,\n)$ if the initial degree of a socle element $f$ does not equal ${\rm socdeg}(A)$, then $f$ must be part of a minimal generating set of $\n$.

\begin{Lemma}\label{Lstretched}$($\cite[3.2]{EV}$)$
Let $(A,\n)$ be a stretched Artinian  local ring with socle degree $s$ and let $f\in\n$ be a socle element. If $f\notin \n^s$ then $f\notin \n^2$.
\end{Lemma}

\demo If $s\,\leq\,2$ we are done. So we may assume that $s\geq 3$.
By assumption of stretchedness, one has the following Hilbert function for $A$:
$$ HF_A\;:\quad 1 \quad c \quad 1 \quad 1 \quad \dots \quad 1 \quad 0_{\longrightarrow}$$
Assume $f\in\n^2$ and set by $i$ the minimun integer with $f\notin\n^{i+1}$. Since $f\notin \n^s$ we have $2\,\leq\,i\,\leq\,s-1$. Then $HF_{A}(i)=1$, which implies  $\n^i=(f)$. Since $f\cdot \n = (0)$, then $\n^{i+1}=\n^i\cdot \n =f\cdot \n =(0)$. Hence $i+1\,\geq\,s+1$, contradicting that $i\,<\,s$. Therefore $f\notin\n^2$.
\QED

\bigskip

A classical result of Sally \cite{Sa} gives the structure of stretched Gorenstein Artinian  local algebras. Recently, Elias and Valla  \cite{EV} generalized it to {\em any} stretched Artinian  local algebras as follows:

\begin{Theorem}\label{SS}$($\cite[1.2]{Sa}, \cite[3.1]{EV}$)$
Let $(R,\m)$ be a $c$-dimensional regular local ring with  ${\rm char}\,R/\m =0$. Let $I\subseteq \m^2$ be an $\m$-primary ideal with $R/I$ stretched of socle degree $s$. Write $\tau(R/I)=r+1$ for some $0\leq r\leq c-1$. Then there exist minimal generators $x_1,\dots,x_c\,$ for the maximal ideal $\m$ and elements $u_{r+1},\dots,u_{c-1} \notin \m$ with
$$I=(x_1\m,\ldots,x_r\m)+J$$
where \begin{equation*}
J = \left\{
\begin{array}{ll}
\!\!(x_{r+i}x_{r+j}\,|\, 1\leq i< j \leq  c-r) + ( x_c^s - u_{r+i}x_{r+i}^2\,|\,1\leq i\leq c-r-1), & \text{if }\;\, r<c-1;\\
\!\!(x_c^{s+1}), & \text{if }\;\, r=c-1.\\
 \end{array} \right.
\end{equation*}
%$J=(x_{r+i}x_{r+j}\,|\, 1\leq i< j \leq  c-r) + ( x_c^s - u_{r+i}x_{r+i}^2\,|\,1\leq i\leq c-r-1)$, if $r< c-1$ and $J=(x_c^{s+1})$, if $r=c-1.$
\end{Theorem}

\medskip

%In the homogeneous case the situation is much easier, since up to graded isomorphism there is only one possible stretched ideal of socle degree $s$. In the case of characteristic zero, it follows by the above theorem, however we give below a short proof that does not require such assumption.
%\begin{Proposition}\label{stretchedh}
%Let $S=k[X_1,\dots,X_c]$ be a polynomial ring over a field. Let $\m=(X_1,\dots,X_c)$ and $I$ a homogeneous $\m$-primary ideal with $R/I$ stretched and socle degree $s\geq 3$. After (possibly) a homogeneous change of coordinates, $I=(X_1\m,\dots,X_{c-1}\m,X_c^{s+1})$.
%\end{Proposition}

%\demo
%Since $s\,\geq\,3$ we have $HF_{S/I}(2)=HF_{S/I}(3)=1$ while $HF_{S/I}(1)=(c-1)+1$. By a Theorem of Cho-Iarrobino (see for instance \cite[Theorem~3.4]{GHMS})
 %this forces $R/I$ to have $c-1$ homogeneous socle elements of degree $1$. After a homogeneous change of coordinates we may assume they are $X_1,\dots,X_{c-1}$. Notice that $X-c^{s+1} \in I$ since $\m^{s+1}\subseteq I$, so $(X_1\cdot\m,\dots, X_{c-1}\cdot m,X_c^{s+1}\subseteq I$ and since these two ideals have the same Hilbert Function they must be equal.
%\QED
%\medskip

In a forthcoming paper we are actually able to slightly improve \ref{SS}, showing that the structure theorem holds for any regular local ring with ${\rm char}\,R/\m \neq 2$. Hence, this will be our only requirement on the residue field in the following.

In a Noetherian local ring $(R,\m)$, let $x_1,\dots,x_r$ be a regular sequence in $\m$; we recall that a {\em monomial} in $x_1,\dots,x_r$ is an element in $R$ of the shape $x_1^{a_1}x_2^{a_2}\cdots x_r^{a_r}$. An $R$-ideal $I$ is dubbed {\em monomial} (in $x_1,\dots,x_r$) if it can be generated by monomials in $x_1,\dots,x_r$.

As an application of Theorem \ref{SS} we prove that the square of a stretched $\m$-primary  ideal $I$ can be very tightly estimated by using a monomial ideal. This fact will be crucial as it will allow us to compute explicitly, or estimate tightly, $e(R/I^2)$ (see the comment after \ref{Tstretched}, or the proof of \ref{c+4}).

\begin{Proposition}\label{Isquare}
Let $(R,\m)$ be a $c$-dimensional regular local ring with ${\rm char}\,R/\m \neq 2$. Let $I\subseteq \m^2$
be an $\m$-primary ideal with $R/I$ stretched of socle degree $s$.
%of socle degree $s\geq 3$.
Write $\tau(R/I)=r+1$ for some  $0\leq r\leq c-1$.
 Then there exist minimal generators $x_1,\dots,x_c\,$ for the maximal ideal $\m$ such that
$$I^2\subseteq\,L=(x_1,\dots,x_{c-1})H + (x_1,\dots,x_{c-1}) x_c^{s+1} + (x_c^{2s}),$$
where $H$ is the monomial ideal  generated by  all  monomials of degree $3$ in $x_1,\dots,x_c$ except for $x_c^3$.
The above inclusion is strict if and only if $r\geq c-2$, and in this case $\lambda(R/I^2)\geq \lambda(R/L)+2$.
%(c) If $r=c-2$ then after choosing a proper minimal generating set $x_1,\dots,x_c$ of $\m$, one has $I^2=J + (x_{c-1}x_c^{s+1}-u_{c-1}x_{c-1}^3x_c,x_c^{2s}+u_{c-1}^2x_{c-1}^4)$ where $J$ is as above.
%(d) If $r=c-1$ then $I$ is monomial $I=(x_1\cdot\m,\dots,x_{c-1}\cdot\m,x_c^{s+1})$ and so $I^2$ is too.
\end{Proposition}

\demo
By Theorem \ref{SS}, there exist minimal generators $x_1,\dots,x_c$ for $\m$ and units $u_{r+1},\dots,u_{c-1}$ with $I=(x_1\m,\ldots,x_r\m)+J$, where $J$ is as in  Theorem \ref{SS}. It is now easy to verify that every element in $I^2$ is a linear combination of the elements of $L$, proving the inclusion $I^2\subseteq L$.

To prove the second part of the statement,  observe that if $r\,\leq\,c-3$ then both $x_c^{s}-u_{c-2}x_{c-2}^2$ and $x_c^s-u_{c-1}x_{c-1}^2$ are in $I$, hence their product is in $I^2$. Since $x_{c-2}^2x_{c-1}^2\in I^2$ and $x_{c-i}^2x_c^s = (x_{c-i}x_c)^2\, x_c^{s-2} \in I^2$ for any $i=1,2$, we have $x_c^{2s} \in I^2$ as well. Therefore $$(x_1\m,\dots,x_{r}\m,x_{r+i}x_{r+j}\,|\,1\,\leq\,i\,<j\,\leq c-r)^2 + (x_c^{2s}) \subseteq I^2.$$
Again, for any $1\leq i\leq c-r-1$,  the fact that $(x_c^s-u_{r+i}x_{r+i}^2)^2\in I^2$  implies that $x_{r+i}^4\in I^2$, since both $x_c^{2s}$ and $x_{r+i}^2x_c^s $ are in $I^2$. Now let $1\leq i\leq c-1$, there exists $1\leq j\leq c-r-1$ such that $i\neq r+j$. Since both $x_{i}x_c(x_c^s-u_{r+j}x_{r+j}^2)\in I^2$ and $x_{i}x_cx_{r+j}^2=(x_{i}x_{r+j})(x_{r+j}x_c) \in I^2$, we must have that $x_{i}x_c^{s+1}\in I^2$.  Finally, for any $1\leq i \leq c-r-1$ and $j\neq r+ i$, since $x_{r+i}x_j(x_c^s-u_{r+i}x_{r+i}^2) \in I^2$, we have $x_{r+i}^3x_j \in I^2$. These arguments show that all the generators of $L$ are in $I^2$ thus proving the equality in the case $r\leq c-3$.

Now assume $r=c-2$. Then $L=(M, x_{c-1}^4, x_{c-1}^3x_c,  x_{c-1}x_c^{s+1}, x_c^{2s})$ and $I^2=(M, x_{c-1}x_c^{s+1}-u x_{c-1}^3x_c, x_c^{2s}+u^2 x_{c-1}^4)$, where $u$ is a unit and $M=(x_1,\dots,x_{c-2})H+ (x_1,\dots,x_{c-2}) x_c^{s+1}+ (x_{c-1}^2x_c^2)$ is a monomial ideal.
Observe that $L/I^2=(\overline{x_{c-1}^4}, \overline{x_{c-1}^3x_c})$ and none of the two generators of $L/I^2$ is redundant.
%elements $x_{c-1}^3x_c, x_{c-1}^4, x_{c-1}x_c^{s+1}, x_c^{2s}$  is redundant in the above description of $L$ and none of the elements $x_{c-1}x_c^{s+1}-u x_{c-1}^3x_c, x_c^{2s}+u^2 x_{c-1}^4$ is redundant in the above description of $I^2$.
This shows  that $I^2 \neq L$ and  $\lambda(R/I^2)\geq \lambda(R/L)+2$.
 %where $\mu(J)$ denotes the minimal number of generators of an ideal $J$. The latter inequality in particular yields the inequality of the lengths.

Finally if $r=c-1$ then the claim is easily proved, since $I=(x_1,\dots,x_{c-1}) \m + (x_c^{s+1})$ and one can easily compare the two monomial ideals $I^2$ and $L$.
%$I^2=J + (x_{c-1}x_c^{s+1}-u_{c-1}x_{c-1}^3x_c,x_c^{2s}+u_{c-1}^2x_{c-1}^4)$ where $J$ is as above.
%The proof of the case $r\leq c-3$ is left as a straightforward (but possibly tedious) exercise for the reader.
%(a) is clearly a consequence of (b), (c) and (d), which are left as straightforward (but possibly tedious) exercises for the reader.
\QED

\bigskip

%\noindent We show an example of why we do need to assume that $r+1\neq c-1$.
%\begin{Example}
%Assume $c=3$ and $r=1$. Then by Theorem \ref{SS}, $I=(x_1\cdot\m,x_2x_3,x_3^s-u_2x_2^2)$. Then $I^2=(x_1^2\cdot\m^2,x_1x_2x_3\cdot\m,x_1^2x_3^s,x_1x_3^{s+1},x_1x_2^3,x_2^2x_3^2,x_2x_3{s+1}-u_2x_2^3x_3,x_3^{2s}+u_2^2x_2^4)$ which is clearly not a monomial ideal in $x_1,x_2,x_3$.
%\end{Example}

%Assume $c=4$ and $r=1$. Then by Theorem \ref{SS}, $I=(x_1^2,x_1x_2,x_1x_3,x_1x_4,x_2x_3,x_2x_4,x_3x_4,x_4^s-u_2x_2^2,x_4^s-u_3x_3^2)$.
%By Proposition \ref{Isquare} $I^2=(x_1,x_2,x_3)(x_1^3,x_1^2,x_1^2x_3,x_1^2x_4,x_1x_2^2,x_1x_2x_3,x_1x_2x_4,x_1x_3^2,x_1x_3x_4,x_1x_4^2,x_2^3,x_2^2x_3,x_2^2x_4x_2x_3^2,x_3^{s+1}) + (x+4^{2s})$
%Assume $c=3$ and $r=1$. Then by Theorem \ref{SS}, $I=(x_1^2,x_1x_2,x_1x_3,x_2x_3,x_3^s-u_2x_2^2)$.
%By Proposition \ref{Isquare} $I^2=(x_1,x_2)(x_1^3,x_1^2x_2,x_1^2x_3,x_1x_2^2,x_1x_2x_3,x_1x_3^2,x_2^3,x_2^2x_3,x_2x_3^2,x_3^{s+1})+(x_3^{2s})$

%Assume $c=4$ and $r=1$. Then by Theorem \ref{SS}, $I=(x_1^2,x_1x_2,x_1x_3,x_1x_4,x_2x_3,x_2x_4,x_3x_4,x_4^s-u_2x_2^2,x_4^s-u_3x_3^2)$.
%By Proposition \ref{Isquare} $I^2=(x_1,x_2,x_3)(x_1^3,x_1^2,x_1^)$

The next theorem is the main result of this section.  It shows that unless %the socle degree is $2$ or
the embedding codimension is $3$, the first conormal module (hence square) of any  generically a complete intersection stretched Cohen-Macaulay ideal is not Cohen-Macaulay.

\begin{Theorem}\label{Tstretched}
Let $(R,\m)$ be a regular local ring with  ${\rm char}\,R/\m \neq 2$ and  $I$ a stretched Cohen-Macaulay ideal
that is generically a complete intersection. % but not a complete intersection. is  %with socle degree $s\geq 3$
Assume either $R$ is containing a field or $I\subseteq \m^2$ and write $c={\rm ecodim}(R/I)$.

$($a$)$  If $c=3$ and $I/I^2$ $($or $R/I^2$$)$ is Cohen-Macaulay then $R/I$ is Gorenstein.

$($b$)$ If $c\geq 4$ then $I/I^2$ and $R/I^2$ are not Cohen-Macaulay.
\end{Theorem}

\demo
Notice that if $R$ contains a field, by Proposition \ref{NonDeg2} we can reduce at once to the case where $I\subseteq \m^2$. Hence, in either case we may assume that $I\subseteq \m^2$ and then $c={\rm ht}\,I$.
Let $s={\rm socdeg}(R/I)$. Since $R/I$ is stretched, we have that $e(R/I)=c+s$. Hence, by Proposition \ref{Generalized},\, $e(R/I^2)=(c+1)(c+s)$.
We will show by contradiction  that $I/I^2$ is not Cohen-Macaulay if $c\geq 4$ or if $c=3$ and $R/I$ is not Gorenstein. So assume   $I/I^2$ is  Cohen-Macaulay. Since $R/I$ is  also Cohen-Macaulay, the short exact sequence
$$ 0 \longrightarrow I/I^2 \longrightarrow R/I^2 \longrightarrow R/I \longrightarrow 0$$
shows that $R/I^2$ is  Cohen-Macaulay as well. Hence, after possibly extending the residue field of $R$, there exists an $R$-ideal $J$ whose images in $R/I$ and $R/I^2$ are minimal reductions of $\m_{R/I}$ and $\m_{R/I^2}$, respectively, generated by systems of parameters. Let $\overline{R}$ be the regular local ring $\overline{R}=R/J$.
Then %$\lambda(R/(I, \underline y))=e(R/I)$ and
$\lambda(\overline{R}/\overline{I}^2)=\lambda(R/I^2\!+\!J)=e(R/I^2)=(c+1)(c+s)$. %In order to find a contradiction we compute $\lambda(R/(I^2, \underline y))$ in a different way to show that we get a different result.
However, by assumption  $\overline{R}/\overline{I}$ is a stretched Artinian algebra of socle degree $s$, hence by Proposition \ref{Isquare} there exist minimal generators $x_1,\dots, x_c$ of $\overline{\m}$ such that
%$(I, \underline y)=( x_1\m,\ldots,x_r\m,
% J, \underline{y})$, where $r\leq c-1$ and $J=(x_{r+i}x_{r+j}\,|\, 1\leq i< j \leq  c-r) + ( x_c^s - u_{r+i}x_{r+i}^2\,|\,1\leq i\leq c-r-1)$.\\
%By Proposition \ref{Isquare},
$\overline{I}^2\subseteq L$, where  $L$ is the $\overline{R}$-ideal defined as $L=(x_1,\dots,x_{c-1})H + (x_1,\dots,x_{c-1}) x_c^{s+1} + (x_c^{2s})$, and  $H$ is the monomial ideal in $\overline{R}$ generated by  all the monomials of degree $3$ in $x_1,\dots,x_c$ except for $x_c^3$. Now, it is easy to see that the monomials $x_ix_c^{t}$ and $x_c^{j}$ are not in $L$ for any $1\leq i \leq c-1$, $3\leq t\leq s$ and  $4\leq j\leq 2s-1$.
Since $L$ is a monomial ideal in $x_1,\dots,x_c$,
the initial forms of the monomials $x_ix_c^{t}$ and $x_c^{j}$, where $1 \leq i \leq c-1$, $3\leq t\leq s$ and  $4\leq j\leq 2s-1$,
%are not in $\overline{L}$, then they
are actually $R/\m$-linearly independent in ${\rm gr}_{\overline{\m}}(\overline{R}/L)$ (see also the comment at the end of the proof).
Therefore, if one writes
$$1 \qquad c \qquad h_2 \qquad h_3 \qquad \dots$$
for the Hilbert function of $\overline{R}/ L$, the above statement shows that $h_{j+1}\geq c$ for every $3\leq j\leq s$ and $h_{j}\geq 1$ for every $s+2 \leq j\leq 2s-1$.% (see also the comment at the end of this proof).
Then $\lambda(\overline{R}/L)\geq 1 + c + \binom{c+1}{2} + \binom{c+2}{3} + c(s-2) + s-2$. By computations, if $c\geq 4$ then $1 + c + \binom{c+1}{2} + \binom{c+2}{3} + c(s-2) + s-2> (c+1)(c+s)$. Hence $\lambda(\overline{R}/\overline{I}^2)\geq \lambda(\overline{R}/L)> (c+1)(c+s)$ that is a contradiction. Assume then $c=3$. Since $R/I$ is not Gorenstein, then $\overline{R}/\overline{I}$ is not Gorenstein as well and we have $c-2=1\leq r=\tau(\overline{R}/\overline{I})-1$. Now by Proposition \ref{Isquare}, we  have $\lambda(\overline{R}/\overline{I}^2)>\lambda(\overline{R}/L)\geq (c+1)(c+s)$, giving the desired contradiction.
\QED
\bigskip

In the proof, the key point is to estimate the Hilbert function of $\overline{R}/\overline{I}^2$ using the existence of some $R/\m$-linearly independent elements in $\overline{R}/L$. We would like to explain this step better.

In the local case it is usually hard to estimate the Hilbert function at a given degree using linearly independent elements. In fact, even if there are $A/\n$-linearly independent  elements $a_1,\dots, a_n$ of $\n^i \backslash \n^{i+1}$ in an Artinian local ring $(A,\n)$, their initial forms in the associated graded ring  $A^*={\rm gr}_{\n}(A)$ need not be $A/\n$-linearly independent. %Indeed, it might happen that $a_j - a_k \in \m^{i+1}$ for some $j\neq k$.
For instance, the elements $x+y^2,\, x+xy,\,x+x^2$ are clearly $k$-linearly independent in $A=k[x,y]_{(x,y)}/(x,y)^3$,   but  their initial forms in $A^*$ are all equal!
However, let $(R,\m)$ be a regular local ring and fix a regular system of parameters $\underline{x}=x_1,\dots,x_c$, so that every `monomial' item in the following is referred to the regular sequence $\underline{x}$. If $L$ is an $\m$-primary monomial ideal and there are distinct monomial elements $a_1,\dots, a_n$ in $R$ with $a_j \notin L$ for every $j$ and all $a_j$'s have the same initial degree $i$, then their initial forms %in the associated graded ring  $R^*/I^*$
 are $R/\m$-linearly independent, proving that $HF_{R/L}(i)\geq n$.
 In the proof of Theorem \ref{Tstretched},  we use that $\overline{I}^2$ is contained in the monomial ideal $L$ to translate the problem into estimating the Hilbert function of $L$, which is a much easier task by the above. This approach is successful because $L$ estimates very tightly $\overline{I}^2$, as explained by Proposition \ref{Isquare}.

As a consequence of Theorem \ref{Tstretched} and Lemma \ref{LC},
we can now give a positive answer to Question \ref{Vasc2} for stretched algebras.

\begin{Corollary}\label{Vstretched}
%Let $(R,\m)$ be either a regular local ring with
Assume ${\rm char} R/\m\! \neq\! 2$. Then Question \ref{Vasc2} is true for $I$ defining a stretched algebra.
\end{Corollary}
%\demo By Lemma \ref{LC}, we can assume that $c={\rm ecodim}(R/I)\geq 3$. The statement just follows from Theorem \ref{Tstretched}. \QED
%Set, as usual, $s=s(R/I)$. If $s\geq 3$ we are done by Theorem \ref{Tstretched}. So we may  assume that $s=2$ and by Proposition \ref{NonDeg2} we can further assume $I\subseteq \m^2$. By Proposition \ref{Generalized}, $e(R/I^2)=(c+1)e(R/I)=(c+1)(c+2)=c^2+3c+2$. Since $I/I^2$ is Cohen-Macaulay, as in the proof of Theorem \ref{Tstretched},
%we can find a regular sequence $\underline y$ in $R$ generating a minimal reduction of $\m_{R/I}$ and $\m_{R/I^2}$.
%Then %$\lambda(R/(I, \underline y))=e(R/I)=c+2$ and
%$\lambda(R/(I^2, \underline y))=e(R/I^2)=c^2+3c+2$. We show that if $R/I$ is not Gorenstein then $\lambda(R/(I^2, \underline y))>c^2+3c+2$, giving a contradiction. To prove it, write $\m=(x_1, \ldots, x_c, \underline{y})$. If $c\geq 4$ then $\lambda(R/(I^2, \underline y))\geq \lambda(R/((x_1, \ldots, x_c)^4,\underline y))=\frac{c^3}{6} + c^2 + \frac{11}{6}c + 1$. This number is strictly bigger than $c^2+3c+2$ if $c\geq 4$. For $c=3$, however, $\lambda(R/((x_1, x_2, x_3)^4,\underline y))=20=3^3+3\cdot3 + 2$, hence in this latter case it is enough to prove that $(I^2,\underline y)\neq ((x_1, \ldots, x_c)^4,\underline y)$.
%Now Lemma \ref{socdeg2} finishes the proof.
%\medskip
%Notice, however, that Theorem \ref{Tstretched} is actually stronger than Question \ref{Vasc2}, in the case of stretched algebras.
\medskip

\section{Monomial ideals}
In the present section we prove that Question \ref{Vasc2} has a positive answer for monomial ideals.

Recall that in a regular local ring $(R,\m)$, fixed a regular system of parameters $\underline x=x_1,\dots,x_c$, an element $a\in R$ is said to be a {\em monomial} in $\underline x$, if there exist non negative integers $a_1,\dots,a_c$ with $a=x_1^{a_1}x_2^{a_2}\cdots x_c^{a_c}$. An $R$-ideal $I$ is said to be a {\em monomial ideal} in $\underline x$, if it can be generated by monomial elements in $\underline x$.

Rinaldo, Terai and Yoshida recently proved the following result.

\begin{Theorem}$($\cite{RTY}$)$\label{RTY}
Let $\Delta$ be a simplicial complex on $V = [n]$, and $I_{\Delta} \subseteq  S = K[X_1, . . . , X_n]$ be
the Stanley-Reisner ideal of $\Delta$. If $S/I_{\Delta}^2$ is Cohen-Macaulay for any field $K$, then $I_{\Delta}$ is Gorenstein.\\ Equivalently, 
Question \ref{Vasc2} holds true for any squarefree monomial ideal $J$ with the property that $S/J^2$ is Cohen-Macaulay for any field $K$.
\end{Theorem}

%\begin{Theorem}$($\cite[2.2]{RTY}$)$\label{RTY}
%Let $\Delta$ be a simplicial complex on V = [n], and let $I_{\Delta} \subseteq  S = K[X_1, . . . , X_n]$ denote
%the Stanley-Reisner ideal of $\Delta$. If $S/I_{\Delta}^2$ is Cohen-Macaulay for any field $K$, then $I_{\Delta}$ is Gorenstein.
%\end{Theorem}
%Notice that obviously, one can replace $S$ by $\hat{S}$ and still the same statement will be true.

We can now prove that Question \ref{Vasc2} holds true for {\em any} monomial ideal.
% and, at the same time, we generalize the result of Rinaldo, Terai and Yoshida.
%Given a power series ring over a field $k[\![\underline x]\!]$, given a monomial ideal $J$, we
%Given a regular local ring $(R,\m)$, containing a field $k$, fix a regular system of parameters $\underline x$. Notice that $\hat{R}=k[\![\underline x]\!]$. We define the class $\mathbb H_{\underline x}$ of ideals $J$ in $R$ with the property that $J$ is a monomial ideal in $\underline x$ and its image in $K[\![\underline x]\!]$ is Cohen-Macaulay for any field $K$.
\begin{Theorem}\label{monomial}
In the setting of Question \ref{Vasc2}, further assume that $R$ contains a field $k$ and $I$ is a monomial ideal in a regular system of parameters $\underline x$ whose corresponding monomial ideal $\tilde{I}$ in $\tilde{R}=K[\![\underline x]\!]$ is Cohen-Macaulay for any field $K$. If either $\tilde{I}/\tilde{I}^2$ or $\tilde{R}/\tilde{I}^2$ is Cohen-Macaulay for any field $K$ then $R/I$ is Gorenstein.%. Let $\underline x$ be a regular system of parameters of $R$ and $I$ be a monomial ideal whose corresponding monomial ideal $\tilde{I}$ in $\tilde{R}=K[\![\underline x]\!]$ is Cohen-Macaulay for any field $K$. If $\tilde{I}/\tilde{I}^2$ (equivalently, $\tilde{R}/\tilde{I}^2$) is Cohen-Macaulay for any field $K$ then $R/I$ is Gorenstein.
%
%In the setting of Question \ref{Vasc2}, assume $R$ contains a field and $I$ is a monomial ideal in a regular system of parameters $\underline x$ with $I \in \mathbb H_{\underline x}$. If $I^2 \in H_{\underline x}$ then $R/I$ is Gorenstein.
\end{Theorem}

\demo
Without loss of generality we may assume that $R$ is complete, hence we may assume that $R=k[\![\underline x]\!]$ is a power series over a field $k$.
Let $J\subseteq S=k[\![\underline x,\underline y]\!]$ be a squarefree monomial ideal obtained by polarization from $I$. Then, $(S,J)$ is a deformation of $(R,I)$. Fix a field $K$. Set %$\tilde{R}= K[\![\underline x]\!]$,
$\tilde{S}=K[\![\underline x,\underline y]\!]$ and let $\tilde{I}$ and $\tilde{J}$ be the corresponding monomial ideals in $\tilde{R}$ and $\tilde{S}$, respectively. Notice that $\tilde{R}/\tilde{I}$ and $\tilde{R}/\tilde{I}^2$ are Cohen-Macaulay, by our assumption on $I$ and $I^2$
(if $\tilde{I}/\tilde{I}^2$ is Cohen-Macaulay, then by the short exact sequence displayed after the statement of Question \ref{Vasc2} we also have that $\tilde{R}/\tilde{I}^2$ is Cohen-Macaulay).

Since $\tilde{I}$ is generically a complete intersection, a proof similar to the one of \cite[Theorem~2.1]{Hu} shows that $\tilde{S}/\tilde{J}$ and  $\tilde{S}/\tilde{J}^2$ are Cohen-Macaulay for any field $K$, because $\tilde{R}/\tilde{I}$ and $\tilde{R}/\tilde{I}^2$ are.
We have then just proved that $J$ is a squarefree monomial ideal with $\tilde{S}/\tilde{J}^2$ Cohen-Macaualy for any field $K$. Then, by Theorem \ref{RTY} (applied to $R=k[\![\underline x]\!]$), $S/J$ is Gorenstein. Finally, since $(S,J)$ is a deformation of $(R,I)$ and $S/J$ is Gorenstein then $R/I$ is Gorenstein, finishing the proof.
\QED
\bigskip
\begin{Corollary}
With notation as in \ref{monomial}, Question \ref{Vasc2} holds true for all monomial ideals $I$ with the property that $\tilde{R}/\tilde{I}$ and $\tilde{I}/\tilde{I}^2$ are Cohen-Macaulay for any field $K$.

\end{Corollary}
%\begin{Corollary}
%Question \ref{Vasc2} is true if $R$ contains a field $k$ and $I$ is a monomial ideal in a regular system of parameters $\underline x$ such that
%$\tilde{I}$ and $\tilde/\tilde{I
% corresponding monomial ideal $\tilde{I}$ in $\tilde{R}=K[\![\underline x]\!]$ is Cohen-Macaulay for any field $K$. If either $\tilde{I}/\tilde{I}^2$ or $\tilde{R}/\tilde{I}^2$ is Cohen-Macaulay for any field $K$ then $R/I$ is Gorenstein.
%. Let $\underline x$ be a regular system of parameters of $R$ and $I$ be a monomial ideal whose corresponding monomial ideal $\tilde{I}$ in $\tilde{R}=K[\![\underline x]\!]$ is Cohen-Macaulay for any field $K$. If $\tilde{I}/\tilde{I}^2$ (equivalently, $\tilde{R}/\tilde{I}^2$) is Cohen-Macaulay for any field $K$ then $R/I$ is Gorenstein.
%\end{Corollary}
%
%An immediate consequence is that Question \ref{Vasc2} has a positive answer for monomial ideals and generalizes the result of Rinaldo, Terai and Yoshida.
%\begin{Corollary}
%In the setting of Question \ref{Vasc2}, further assume that $R$ contains a field $k$. Let $\underline x$ be a regular system of parameters of $R$ and $I$ be a monomial ideal whose corresponding monomial ideal $\tilde{I}$ in $\tilde{R}=K[\![\underline x]\!]$ is Cohen-Macaulay for any field $K$. If $\tilde{I}/\tilde{I}^2$ (equivalently, $\tilde{R}/\tilde{I}^2$) is Cohen-Macaulay for any field $K$ then $R/I$ is Gorenstein.
%\end{Corollary}

\section{Short algebras}

In this section we will study the  conormal module of ideals defining short algebras. We will prove that  for these ideals the  conormal module (equivalently, the square) is almost never Cohen-Macaulay. Again, the information on the Hilbert function of such algebras will play a crucial role.

Recall that in a  Noetherian local ring $(R,\m)$, the {\em associated graded ring} of $R$ with respect to $\m$ is defined as the graded algebra $R^*={\rm gr}_{\m}(R)=\bigoplus_{i\,\geq\,0}\m^i/\m^{i+1}$.
 For an element $0\neq a\in R$, there is a unique integer $i$ with $a \in \m^i\backslash\m^{i+1}$. This integer is called the {\em initial degree} of $a$, denoted ${\rm indeg}(a)$. The {\em initial form} of $a$ is $a^*=(a+\m^{i+1})/\m^{i+1} \in R^*$. Similarly, the {\em ideal of initial forms} of  an $R$-ideal $J$ is defined as  $J^*=( a^* \,|\, a\in J)R^*$. It is well known that ${\rm gr}_{\m}(R/J)\simeq R^*/J^*$.
The definition of short algebra can be found, for instance, in \cite{CRV}.

\begin{Definition}\label{defshort}
An Artinian local algebra $(A,\n)$ is {\rm short} if there exist integers $c$ and $s$ %$c\in\mathbb Z_+$ and $s\in\mathbb Z_+$
with $HF_A(j)=HF_{A/\n[X_1,\dots,X_c]}(j)$ for every $j<s$, and $HF_A(s+1)=0$ $($i.e., if $A^*$ is a short graded algebra$)$.
A Cohen-Macaulay local ring $(R,\m)$ is {\rm short} if $R/J$ is short for any   minimal reduction $J$ of $\m$.
\end{Definition}

The second part of the next remark will be used to deal with Artinian short algebras, as it shows that one can reduce any computation on the length of $R/I$ or $R/I^2$ to the homogeneous case.
%Let $(A,\n)$ be a Artinian local ring, following \cite{Em} we say that $A$ is {\it canonically graded} if $A$ is isomorphic to the completion (with respect to the homogeneous maximal ideal) of its associated graded ring, i.e., if $A\simeq \widehat{{\rm gr}_{\n}(A)}$.

\begin{Remark}\label{homog}
Let $(R,\m, k)$ be a $d$-dimensional regular local ring containing a field and  $I$ an $\m$-primary ideal. Then $R/I\simeq \tilde{R}/\tilde{I}$, where $\tilde{R}=k[x_1,\dots,x_d]$ and $\tilde{I}$ is generated by polynomials of degree at most $s+1$, where $s={\rm socdeg}(R/I)$.
Furthermore, if $I$ defines a short algebra, %Then for some positive integer $s$, $\m^{s+1}\subseteq I \subseteq \m^s$. Hence $R/I$ is canonically graded, in particular, $(I^*)^n=(I^n)^*$ for any $n\,\geq\,1$ and  $\lambda(R/I^n)=\lambda(R^*/(I^*)^n)$.
%$\lambda(R/I^2)\,>\,(c+1)\lambda(R/I)$ if and only if $\lambda(R^*/(I^*)^2)\,>\,(c+1)\lambda(R^*/I^*)$.
then $\tilde{I}$ can be taken to be generated by homogeneous polymonials.  In particular, $(I^*)^n=(I^n)^*$ for any $n\,\geq\,1$ and  $\lambda(R/I^n)=\lambda(R^*/(I^*)^n)$. %Finally, if $I\subseteq \m^2$ then $\m^{s+1}\subseteq I \subseteq \m^s$, where $s={\rm socdeg}(R/I)$.
\end{Remark}

\demo
For the first part of the statement, since $R/I$ is Artinian, then $R/I\simeq \widehat{R/I}\simeq \hat{R}/\hat{I}$. Since $R$ contains a field, one has $\hat{R}\simeq k[\![x_1,\dots,x_d]\!]$, and since $I$ is $\m$-primary with $\m^{s+1}\subseteq I$, one may assume that $\hat{I}=(f_1,\dots,f_r, \m^{s+1})$, where   $f_1,\dots,f_r$ are of degrees at most $s$. However, since $k[x_1,\dots,x_d]/(f_1,\dots,f_r, \m^{s+1})$ is Artinian, one has that $k[x_1,\dots,x_d]/(f_1,\dots,f_r, \m^{s+1}) \simeq k[\![x_1,\dots,x_d]\!]/(f_1,\dots,f_r, \m^{s+1})$, finishing the proof.

The second part follows from the first one and the definition of short algebra.
\QED
\bigskip

%\begin{Remark}\label{homog}
%Let $(R,\m, k)$ be a regular local ring containing a field and  $I$  an $\m$-primary ideal defining a short algebra. %Then for some positive integer $s$, $\m^{s+1}\subseteq I \subseteq \m^s$. Hence $R/I$ is canonically graded, in particular, $(I^*)^n=(I^n)^*$ for any $n\,\geq\,1$ and  $\lambda(R/I^n)=\lambda(R^*/(I^*)^n)$.
%$\lambda(R/I^2)\,>\,(c+1)\lambda(R/I)$ if and only if $\lambda(R^*/(I^*)^2)\,>\,(c+1)\lambda(R^*/I^*)$.
%Then $R/I\simeq \widehat{R/I}\simeq \hat{R}/\hat{I}$ where $\hat{R}$ is a power series ring over $k$ and $\hat{I}$ is generated by homogeneous polymonials of $\hat{R}$.  In particular, $(I^*)^n=(I^n)^*$ for any $n\,\geq\,1$ and  $\lambda(R/I^n)=\lambda(R^*/(I^*)^n)$.\\
%Furthermore, if $I\subseteq \m^2$, then $\m^{s+1}\subseteq I \subseteq \m^s$, where $s={\rm socdeg}(R/I)$.
%\end{Remark}

For instance, let $(R,\m)$ be a regular local ring of dimension $4$ with $\m=(x_1,\dots,x_4)$. Assume $I=(x_1x_2^3,x_2x_3^2x_4,x_1^3x_2,x_4^4,\m^5)$. Then $R/I$ is a short algebra, with
$$k[X_1,\dots,X_4]/(X_1X_2^3,X_2X_3^2X_4,X_1^3X_2,X_4^4,M^5)\simeq R/I,$$
where $M=(X_1,\dots,X_4)$.

\smallskip

From now on, we will assume that $R$ contains a field.
The following  theorem is the main result of this section. It shows that the conormal module (equivalently, the square) of a short algebra is not Cohen-Macaulay if the socle degree is at least 3. Somewhat surprisingly, in \ref{count} we exhibit examples showing that, if the socle degree is 2, the conormal module (equivalently, the square) can actually be Cohen-Macaulay. These examples give a negative answer to Question \ref{Vasc2}.

\begin{Theorem}\label{short}
Let $R$ be a regular local ring and $I$ a Cohen-Macaulay ideal that is generically a complete intersection.  Let $c={\rm ecodim}(R/I)\geq 2$ and $s={\rm socdeg}(R/I)\geq 3$. If $R/I$ is short then $I/I^2$ $($equivalently, $R/I^2$$)$ is not Cohen-Macaulay.
\end{Theorem}

\demo
As in Theorem \ref{Tstretched}, one can always assume that $I\subseteq \m^2$. Let $J$ be an $R$-ideal, whose image in $R/I$ is a minimal reduction of $\m_{R/I}$, generated by a system of parameters (after possibly extending the residue field of $R$, the existence of $J$ is guaranteed). Now, write $\m=(x_1,\dots,x_c)+J$ for the maximal ideal of $R$, where $x_1,\dots,x_c$ is part of a regular system of parameters of $R$. Let $n_{i}$ be the number of monomials of degree $i$ in a polynomial ring in $c$ variables and  $N_{i}$ the number of monomials of degree $i$ in $c+1$ variables. By elementary combinatorics (or by the theory of Hilbert functions), it is easily checked that
$$\sum_{i=0}^{s-1}n_i=N_{s-1} \mbox{ for any } s \geq 1.$$
Hence $e(R/I)=\sum_{i=0}^{s-1}n_i+\binom{c+s-1}{s}\!-\! q=\sum_{i=0}^{s}n_i-\! q=N_{s} \!-\! q$, where $0\leq q\,<\,\binom{c+s-1}{s}$ is the number of minimal generators of $I$ lying in $\m^s \backslash \m^{s+1}$. To show that $I/I^2$ is not Cohen-Macaulay, it is enough to show that
$$e(R/I^2)=\lambda(R/I^2\!+\!J)>(c+1)(N_{s}\!-\! q).$$
Write $\n=(x_1,\dots,x_c)$. Since $R/I\!+\!J$ is short, then $I+J=(f_1,\dots,f_q,\,\n^{s+1})+J$ where the $f_i$'s are all homogeneous of degree $s$ in the variables $x_1,\dots,x_c$ in the sense of Remark \ref{homog}. Therefore $I^2+J\subseteq \n^{2s}+J$ and
$$\lambda(R/I^2 + J)\geq \lambda(R/\n^{2s}+J)=\sum_{i=0}^{2s-1}n_i=N_{2s-1}.  $$
Hence, to finish the proof, it is enough to show that $N_{2s-1}>(c+1)(N_{s}\!-\! q)$ for all $c\geq 2$ and $s\geq 3$. To do this,  we write
$\frac{N_{2s-1}}{c+1}=\frac{\binom{c+2s-1}{2s-1}}{c+1}=\frac{\Pi_{j=0}^{c-1}\left(\frac{2s+j}{j+1}\right)}{c+1}=\Pi_{j=0}^{c-1}\left(\frac{2s+j}{j+2}\right)$ and  similarly, $N_s=\Pi_{j=0}^{c-1}\left(\frac{s+j+1}{j+1}\right)$.
Then we are down to show that $\Pi_{j=0}^{c-1}\left(\frac{2s+j}{j+2}\right)-\Pi_{j=0}^{c-1}\left(\frac{s+j+1}{j+1}\right)\!+\! q > 0$. Set  $Q(c,s)=\Pi_{j=0}^{c-1}\left(\frac{2s+j}{j+2}\right)-\Pi_{j=0}^{c-1}\left(\frac{s+j+1}{j+1}\right)$, we will show that  $Q(c,s)+q\,>\,0$.
 %where The idea is to show that for almost every $c$ and $s$ $Q(c,s)\,>\,0$, which finishes the proof in such cases. The few remaining cases will be proved by arguing on $A$.\\
%Since $M$ is always a strictly positive integer, to finish the proof it is enough to show that $Q\geq 0$.
%Consider the expression $Q(c,s)=\Pi_{j=0}^{c-1}\left(\frac{2s+j}{j+2}\right)-\Pi_{j=0}^{c-1}\left(\frac{s+j+1}{j+1}\right)$.\\
\medskip

Claim: if $Q(\overline{c},\overline{s})> 0$ for some positive integers $\overline{c}$ and $\overline{s}$, then $Q(c,s)> 0$ for every $c\geq \overline{c}$ and $s\geq \overline{s}$.

Assume the claim is proved. Since $Q(5,4)=6\,>\,0$, $Q(4,5)=17\,>\,0$,    $Q(3,6)=7>0$ and $Q(2, 8)=\frac{1}{3}>0$, it follows from the claim  that $Q(c,s)>0$ either if  $c\,\geq\,5$, $s\,\geq\,4$,  or if $c\geq 4$, $s\geq 5$, or if $c\geq 3$, $s\,\geq\,6$ or if $c\,\geq\,2$, $s\,\geq\,8$. We are left to check the following several cases.

If $c=3$, $s=5$ then $Q(3,5)=-1 < 0$. However, if $q\geq 2$ then $\lambda(R/I^2\!+\!J)-(3+1)(N_{5}- q)\geq(3+1)(Q(3,5)+q)> 0$ and we are done. Assume now $q\leq 1$. One has $$\lambda(R/I^2\!+\!J)\geq \lambda(R/\n^{10}+\!J)+[\binom{c-1+10}{10}-q(q+1)/2]\geq \lambda(R/\n^{10}+\!J)+65.$$ Thus $\lambda(R/I^2\!+\!J)-(3+1)(N_{5}- q)\geq (3+1)(Q(3,5)+q)+65\geq -4+65=61>0$.

If $c=3$, $s=4$ then $Q(3,4)=-5$. In this case, if $q\geq 6$ then $\lambda(R/I^2\!+\!J)-(3+1)(N_{4}- q)\geq(3+1)(Q(3,4)+q)> 0$ and we are done. So we may assume $q\leq 5$. Then,
%$\lambda(R/(I^2, \underline y))\geq \lambda(R/(\n^{8}+\underline{y}))+(45-q(q+1)/2)$. Thus
$\lambda(R/I^2\!+\!J)-(3+1)(N_{4}- q)\geq (3+1)(Q(3,4)+q)+(45-q(q+1)/2)=-1/2(q^2-7q-50)>0$ whenever $0\leq q\leq 5.$

The cases $c=3, s=3$, or $c=4, s=3, 4$, or $c=5, s=3$, or $c=2$ and $3\leq s\leq 7$ can be proved similarly.
Thus, the proof is done upon the condition that we prove the above claim. This is accomplished in the next lemma.
\QED
\bigskip

\begin{Lemma}
For any $c \in \mathbb Z_+$ and $s \in \mathbb Z_+$, define $Q(c,s)=\Pi_{j=0}^{c-1}\left(\frac{2s+j}{j+2}\right)-\Pi_{j=0}^{c-1}\left(\frac{s+j+1}{j+1}\right)$. Then\\
(a) $Q(c,s+1)\,\geq\,Q(c,s);$\\
(b) if $c\,\geq2$ and $s\,\geq\,3$ then $Q(c+1,s)\,\geq\,Q(c,s).$\\
In particular, if $Q(\overline{c},\overline{s})\,>\,0$ for some fixed positive integers $\overline{c}\,\geq\,2$ and $\overline{s}\,\geq\,3$ then $Q(c,s)\,>\,0$ for all $c\,\geq\,\overline{c}$ and $s\,\geq\,\overline{s}$.
\end{Lemma}

\demo
 %For any $j\,\geq\,0$ write
%$A_j=\frac{2s+j}{j+2}$ and $B_j=\frac{s+j+1}{j+1}$, then $Q(c,s)=\Pi_{j=0}^{c-1}A_j - \Pi_{j=0}^{c-1}B_j$. Now
For any two positive integers  $c$ and $s$,
\begin{eqnarray*}
Q(c,s+1)&=& \Pi_{j=0}^{c-1}\left(\frac{2s+2+j}{j+2}\right)-\Pi_{j=0}^{c-1}\left(\frac{s+1+j+1}{j+1}\right)\\
&=& \Pi_{j=0}^{c-1}\frac{2s+2+j}{2s+j}\,\frac{2s+j}{j+2} - \Pi_{j=0}^{c-1}\frac{s+1+j+1}{s+j+1}\,\frac{s+j+1}{j+1}\\
&\geq &\Pi_{j=0}^{c-1}\frac{s+1+j+1}{s+j+1}\,\frac{2s+j}{j+2} - \Pi_{j=0}^{c-1}\frac{s+1+j+1}{s+j+1}\,\frac{s+j+1}{j+1}\\
&= & \Pi_{j=0}^{c-1}\frac{s+1+j+1}{s+j+1} \left(\Pi_{j=0}^{c-1}\frac{2s+j}{j+2} - \Pi_{j=0}^{c-1}\frac{s+j+1}{j+1}\right)\\
&\geq & Q(c,s)
\end{eqnarray*}
%\begin{eqnarray*}
%Q(c,s+1)&=& \Pi_{j=0}^{c-1}\left(\frac{2s+2+j}{j+2}\right)-\Pi_{j=0}^{c-1}\left(\frac{s+1+j+1}{j+1}\right)\\
%&=& \Pi_{j=0}^{c-1}\left(A_j + \frac{2}{j+2}\right) - \Pi_{j=0}^{c-1}\left(B_j + \frac{1}{j+1}\right)\\
%&\geq &\Pi_{j=0}^{c-1}A_j - \Pi_{j=0}^{c-1}B_j\\
%&=&Q(c,s)
%\end{eqnarray*}
where the first inequality holds since $\frac{2s+2+j}{2s+j}\geq \frac{s+1+j+1}{s+j+1}$, for all $j\,\geq\,0$.

Similarly to show part (b), write $$Q(c,s)=\overbrace{\Pi_{j=0}^{c-1}\left(\frac{2s+j}{j+2}\right)}^{Q_1}-
\overbrace{\Pi_{j=0}^{c-1}\left(\frac{s+j+1}{j+1}\right)}^{Q_2}.$$ Then
\begin{eqnarray*}
Q(c+1,s)&=&Q_1(c+1,s)-Q_2(c+1,s)\\
&=&
Q_1(c,s)\,\frac{c+2s}{c+2}-
Q_2(c,s)\,\frac{c+1+s}{c+1}\\
&\geq &Q_1(c,s)-Q_2(c,s)\\
&=&Q(c,s)
\end{eqnarray*}
where the inequality holds since $\frac{c+2s}{c+2}\geq \frac{c+s+1}{c+1}\geq 1$, for all $s\geq 3$ and $c\geq 2$.
\QED
\bigskip

The following proposition deals with the case of socle degree 2. Sharp numerical conditions on the multiplicity are provided to ensure that the conormal module is not Cohen-Macaulay.

The idea of the proof is to assume that $I/I^2$ is Cohen-Macaulay, reduce to the case where $I$ is $\m$-primary and produce a contradiction showing that
\begin{equation}\label{EQ} e(R/I^2)> (c+1)e(R/I).
\end{equation}
If the number $q$ of quadric minimal  generators of $I$ is `large' (part (a)), then the multiplicity of $R/I$ is `small', hence the right-hand side in (\ref{EQ}) is `small' and we can prove the above strict inequality. On the other hand, if the number $q$ of quadric minimal generators  of $I$ is `small' (parts (b)-(c)), then $I^2$  has a fairly small amount of minimal generators in degree $4$, hence the multiplicity of $R/I^2$ becomes `large'. In this case, then, the left-hand side of (\ref{EQ}) is `large' and we are again able to prove the above strict inequality.
Finally, notice that if $q$ is `intermediate' there may be no contradiction, as $I/I^2$ may actually be Cohen-Macaulay -- see Example \ref{count}.

We are grateful to A. Conca for bringing our attention to \cite[Theorem~2.4]{COR}, that finishes the proof of \ref{short2}.

\begin{Proposition}\label{short2}
Let $(R,\m)$ be a regular local ring and $I$ a Cohen-Macaulay ideal that is generically a complete intersection with $s={\rm socdeg}(R/I)=2$ and $c={\rm ecodim}(R/I)\geq 3$. Write $e(R/I)=1+c+\binom{c+1}{2} - q$,  for some non negative integer $q<\binom{c+1}{2}$.  Then $I/I^2$ is not Cohen-Macaulay if $q$ satisfies  one of the following numerical conditions$:$\\
 \vspace{0.02in}
 $($a$)$ $q\,>\,\frac{c^2+2c}{3},$\\
\vspace{0.02in}
$($b$)$ $\frac{2c+1-\sqrt{K}}{2}<q<\frac{2c+1+\sqrt{K}}{2}$, where $K=(2c+1)^2+8(N_4-(c+1)N_2)$ and $N_2, N_4$ are defined
as in the proof of Theorem \ref{short}, \\
\vspace{0.02in}
$($c$)$ $\frac{2c+1-2n_3-\sqrt{K}}{2}<q<\frac{2c+1-2n_3+\sqrt{K}}{2}$, where $K=(2n_3-2c-1)^2+8(N_5-(c+1)N_2)$, where $N_2, N_5, n_3$ are defined
as in the proof of Theorem \ref{short}. \\
Furthermore, $I/I^2$ and $R/I^2$ are not Cohen-Macaulay if $c=3$, $R/I$ is not Gorenstein and ${\rm char}\,R/\m\neq 2$, or if  $c=4$.
%$($d$)$\\
%\vspace{0.02in}
%$($b$)$,\\
%\vspace{0.02in}
%$($c$)$ $q\,>\,\frac{c^2+2c}{3},$\\
%\vspace{0.02in}
%$($d$)$ $\frac{2c+1-\sqrt{K}}{2}<q<\frac{2c+1+\sqrt{K}}{2}$, where $K=(2c+1)^2+8\binom{c-1+4}{4}-\frac{8}{3}(c+1)(c^2+2c)$,\\
%$($d$)$ $\binom{q+1}{2}\,\leq\,\binom{c+3}{4}$ and $q \,< \, \frac{(2c+1) + \sqrt{(2c+1)^2 + 4K}}{2},$ where $K=\frac{c^4-2c^3-13c^2-10c}{12}$,\\\\
%then $I/I^2$ is not Cohen-Macaulay.
\end{Proposition}

\demo
By Proposition \ref{NonDeg2}, we may assume that $I\subseteq\m^2$ and then $c={\rm ht}\,I$.
By Proposition \ref{Generalized}, $e(R/I^2)=(c+1)e(R/I)=(c+1)(N_2-q)$.
 %where $N_2$ is defined as in the proof of Theorem \ref{short}.
 As before, let $J$ be an ideal whose image in $R/I$ gives a minimal reduction of $\m_{R/I}$ generated by a system of parameters.
Then $\lambda(R/I^2\!+\!J)\geq(c+1)(N_2-q),$ where equality holds if and only if $I/I^2$ is Cohen-Macaulay. By computation
$$
\begin{array}{l}
 \lambda(R/I^2\!+\!J)-(c+1)(N_2-q)\\
 \geq  N_3-(c+1)N_2+(c+1)q+{\rm max}\{0,\, n_4-q(q+1)/2\}+{\rm max}\{0,\, n_5-q \,n_3\}.
 %= \frac{1}{6}(c+1)(-2c^2-4c+6q)+{\rm max}\{0, \binom{c-1+4}{4}-q(q+1)/2\}+{\rm max}
 %\{0, \binom{c-1+5}{5}-q\binom{c-1+3}{3}\}.
\end{array}
$$

Now, (a) follows from the fact that
$$\lambda(R/I^2\!+\!J)-(c+1)(N_2-q)\geq N_3-(c+1)N_2+(c+1)q=\frac{1}{6}(c+1)(-2c^2-4c+6q)$$
 which is strictly greater than zero if $q\,>\,\frac{c^2+2c}{3}$.

  Since
\begin{eqnarray*}
\lambda(R/I^2\!+\!J)-(c+1)(N_2-q)&\geq &N_3-(c+1)N_2+(c+1)q+n_4-q(q+1)/2\\
&= &  -\frac{1}{2}[q^2-(2c+1)q-2(N_4-(c+1)N_2)],
\end{eqnarray*}
it is easy to see that  $\lambda(R/I^2\!+\!J)-(c+1)(N_2-q)>0$ if $q$ satisfies the numerical condition given in (b).
Part (c) follows from the inequality
\begin{eqnarray*}
\lambda(R/I^2\!+\!J)-(c+1)(N_2-q)&\geq &-\frac{1}{2}[q^2-(2c+1)q-2(N_4-(c+1)N_2)]+n_5-q\,n_3\\
&= &  -\frac{1}{2}[q^2+(2n_3-2c-1)q-2(N_5-(c+1)N_2)].
\end{eqnarray*}

%$$
%\begin{array}{l}
 %\lambda(R/(I^2, \underline y))-(c+1)(N_2-q)\\
 %\geq  N_3-(c+1)N_2+(c+1)q+{\rm max}\{0, \binom{c-1+4}{4}-q(q+1)/2\}+{\rm max}\{0, \binom{c-1+5}{5}-q\binom{c-1+3}{3}\}\\
 %= \frac{1}{6}(c+1)(-2c^2-4c+6q)+{\rm max}\{0, \binom{c-1+4}{4}-q(q+1)/2\}+{\rm max}\{0, \binom{c-1+5}{5}-q\binom{c-1+3}{3}\}.
%\end{array}
%$$

Finally,  it follows from parts (a), (b) and (c) that
 $\lambda(R/I^2\!+\!J)-(c+1)(N_2-q)>0$ if $c=3,\, q\neq 5$ or if $c=4,\, q\neq 8$. If $c=3,\, q=5$,
  we are  done by  Theorem \ref{Tstretched} (a).
  %and if    we are done again by Theorem \ref{c+3}.
  Assume $c=4,\, q=8$. Set $\overline{R}=R/J$. By Remark \ref{homog}, we may assume that $\overline{R}\simeq k[\![X_1,\dots,X_4]\!]$ and $\overline{I}$ is a homogeneous ideal generated by $8$ quadric polynomials. Since $e(R/I^2)=\lambda(\overline{R}/\overline{I}^2)\geq \lambda(\overline{R}/\overline{\m}^4)=35=(c+1)\lambda(\overline{R}/\overline{I})=(c+1)e(R/I)$, if we show that $\overline{I}^2\neq\overline{\m}^4$, then $e(R/I^2)>(c+1)e(R/I)$, that gives a contradiction. However, since $\overline{I}$ is generated by polynomials, it is enough to show that the ideal $L$ in $k[X_1,\dots,X_4]$ generated by the same $8$ quadrics generating $\overline{I}$ has the property that $L^2\neq \m^4$,  where $\m=(X_1,\dots,X_4)$. This is achieved in \cite[Theorem~2.4]{COR}, where it was shown that any homogeneous ideal $L$ in $k[X_1,\dots,X_4]$ generated by $8$ quadrics has the property that $L^2 \neq \m^4$, if ${\rm char}\,k=0$. Similar methods can be actually employed to show that the same result holds if ${\rm char}\,k \neq 2$, finishing the proof.
%Now assume that ${\rm char}k \neq 0$. To finish the proof it is enough to show that if $I=(q_1,\dots,q_8,\m^3)$ is a homogeneous ideal in $R=k[X_1,\dots,X_4]$, where $q_1,\dots,q_8$ are $k$-linearly independent quadrics and $\m=(X_1,\dots,X_4)$, then $I^2 \neq \m^4$. It is clear that one only needs to prove it in the generic case. Let then $M_1=X_1^2, M_2,\dots,M_{10}=X_4^2$ be the monomials of degree $2$ in $X_1,\dots,X_4$, set $\tilde{q_i}=\sum_{j=0}^{10}Y_{ij}M_j \in \tilde{R}=R[Y_ij]_{ij}$ and write $F_1,\dots,F_{36}$ for the elements $\tilde{q_i}\tilde{q_j}$. Let $N_1=X_1^4,\dots,N_{35}=X_4^4$ be all the monomials of degree $4$ in $R$ and write $\tilde{F_k}=\sum_{m=1}^{35}a_{km}N_m$, where $a_{km}$ are quadrics of degree $2$ in $k[Y_{ij}]_{ij}$ whose coefficients are images of integer numbers. Now, the matrix $A=[a_{km}]_{km}$ has size $36 \,\,x\,\,35$ and to prove the statement in the `generic' case it is enough to show that the rank of $A$ is at most $34$.
%In the case of characteristic zero we know that the rank of $A$ is no more than $34$, and going modulo a prime number $p$ does not increase the rank of the matrix $A$, hence the rank of $A$ is at most $34$, finishing the proof.
\QED

\bigskip

By Example \ref{count}, the numerical conditions given in Proposition \ref{short2}  are actually very sharp. When $c=5$, indeed, $I/I^2$ cannot be Cohen-Macaulay if either $q\!>\!\frac{c^2+2c}{3}=11\frac{2}{3}$ (by Proposition \ref{short2} (a)) or  $q<\frac{2c+1+\sqrt{K}}{2}=11$ (by Proposition \ref{short2} (b)). Hence, the only case left out by our numerical criteria  is $q=11$. However, Example \ref{count} has exactly $c=5$ and $q=11$ and gives a negative answer to Question \ref{Vasc2}!
 %and similar counterexamples for some $c\geq 5$
%This shows that the numerical conditions in Proposition \ref{short2} are sharp when $c=5$. Supported by some numerical data, in Question \ref{counterex} we actually suggest that they are sharp for any $c\,\geq\,5$.\\

A combination of Theorem \ref{short}, Proposition \ref{short2} and Lemma \ref{LC} shows that the question has a positive answer for most short algebras.

\begin{Corollary}\label{VVstretched}
%Let $(R,\m)$ be either a regular local ring with
%Assume that ${\rm char} R/\m=0$. Then
 Question \ref{Vasc2} is true for any ideal defining a short algebra with socle degree at least $3$, or socle degree $2$ and multiplicity satisfying one of the numerical conditions in Proposition \ref{short2}.
\end{Corollary}

%\begin{Remark}\label{Migliore}
%We mentioned in the proof that in the above statement if either $c=3$ or $c=4$ one could prove directly or using the computer that $I/I^2$ is not Cohen-Macaulay. If $c=5$ however, the two parts of Lemma \ref{short2} show that if either $A\,<\,11$ or $A\geq 12$ then $I/I^2$ cannot be Cohen-Macaulay, So the only case left out by our numerical criterion is when $A=11$. We will show later that there actually exists a counterexample for $A=11$ and similar counterexamples showing that Lemma \ref{short2} produces sharp estimates.

\section{Algebras having low multiplicity}

For a Cohen-Macaulay  local ring $(R,\m)$, it is well-known (by Abhyankar's inequality) that the Hilbert-Samuel multiplicity of $R$ with respect to $\m$
is at least $c+1$, where $c$ denotes as usual the embedding codimension of $R$. In this section we show that Question \ref{Vasc2} holds true for any ideal $I$ with $e(R/I)\leq c+4$. In Section 5, we then show that this result is sharp by exhibiting an ideal $I$ with $e(R/I)=c+5$ that is a counterexample to Question \ref{Vasc2}. We will call ideals $I$ with $e(R/I)\leq c+4$  `ideals with low multiplicity'.
%As in the previous section, we always assume that either $R$ contains a field or that we are in the non-degenerate case, i.e., $I\subseteq \m^2$.

We begin this section by proving a general statement that will be used in the next results. After fixing the multiplicity of a Cohen-Macaulay algebra $R/I$, we provide a lower bound on the embedding codimension of $R/I$ to ensure that $I/I^2$ is {\em not} Cohen-Macaulay.

\begin{Proposition}\label{ineq}
Let $R$ be a regular local ring and $I$ a Cohen-Macaulay $R$-ideal that is generically a complete intersection. Assume either $R$ contains a field or $I\subseteq \m^2$. Write $e(R/I)=c+t$ for some positive integer $t$, where $c={\rm ecodim}(R/I)$.
 If $c > \frac{1+\sqrt{1+24(t-1)}}{2}$ then $I/I^2$ $($equivalently, $R/I^2$$)$ is not Cohen-Macaulay.
\end{Proposition}

\demo
Let $\m$ be the maximal ideal of $R$. After possibly invoking Proposition  \ref{NonDeg2} (if $R$ contains a field), we may assume $I\subseteq \m^2$ and then $c={\rm ht}\,I$.
Assume  $I/I^2$ is Cohen-Macaulay, then $R/I^2$ is also Cohen-Macaulay. Let $J$ be an $R$-ideal whose image in $R/I$ is generated by a system of parameters of $\m_{R/I}$. By Proposition \ref{Generalized},
$$\lambda(R/I^2\!+\!J)=e(R/I^2)=(c+1)e(R/I)=(c+1)(c+t)=c^2+(t+1)c+t.$$
Notice  that $I^2\subseteq \m^4$, so $\lambda(R/I^2\!+\!J)\geq\lambda(R/\m^4\!+\!J)=1+c+\binom{c+1}{2}+\binom{c+2}{3}=\frac{1}{6}(c^3+11c)+c^2+1$.
However, $\frac{1}{6}(c^3+11c)+c^2+1 >  c^2+(t+1)c+t$\, if\, $\frac{1}{6}(c^3-c(6t-5)-6(t-1)) > 0$, or  equivalently, if
$(c+1)(c^2-c-6(t-1))>0$.
Since $c$ is positive, it is enough to have $c^2-c-6(t-1)>0$, which holds true if $c > \frac{1+\sqrt{1+24(t-1)}}{2}$.
Hence if $c > \frac{1+\sqrt{1+24(t-1)}}{2}$, then $e(R/I^2)=\lambda(R/I^2\!+\!J)\geq\lambda(R/\m^4\!+\!J)\,>\,(c+1)(c+t)$ which gives a contradiction. Therefore $I/I^2$ can not be Cohen-Macaulay.
\QED
\bigskip

We now employ results from Sections 2 and 3 together with Proposition \ref{ineq} to show that Question \ref{Vasc2} has a positive answer for ideals with multiplicity at most $c+ 4$.  We first deal with ideals of multiplicity at most  $c+3$.

\begin{Theorem}\label{c+3}
Let $(R, \m)$ be a regular local ring containing a field of characteristic $\neq 2$ and $I$ a Cohen-Macaulay $R$-ideal that is generically a complete intersection. Assume that $e(R/I)\,\leq\,c+3$, where $c={\rm ecodim}(R/I)$. If $I/I^2$ $($or $R/I^2$$)$ is Cohen-Macaulay then $R/I$ is Gorenstein.
\end{Theorem}

\demo %Let $\m$ be the maximal ideal of $R$.
Thanks to Lemma \ref{LC} we can assume that $c\geq 3$ and, thanks to Proposition \ref{NonDeg2}, we can always assume $I\subseteq \m^2$ so that $c={\rm ht}\,I$. By Proposition \ref{Generalized},\, $e(R/I^2)=(c+1)e(R/I)$.
Now, let $J$ be an ideal whose image in $R/I$ is a minimal reduction of $\m_{R/I}$ generated by a system of parameters (after extending the residue field such a $J$ always exists). After going modulo $J$, we may further assume that $I$ is $\m$-primary and $\lambda(R/I^2)=e(R/I^2)=(c+1)e(R/I)=(c+1)\lambda(R/I)$. Notice that $I$ may not be generically a complete intersection anymore, but we already used this assumption to estabilish that $e(R/I^2)=(c+1)e(R/I)$ and we will not need it for the rest of the proof.

We now show that if $R/I$ is not Gorenstein then $\lambda(R/I^2)\,>\,(c+1)\lambda(R/I)$,  giving a contradiction. The first case, $e=c+1$, just follows from Proposition \ref{ineq} with $t=1$, because $c\geq 3$. If $e=c+2$, then $R/I$ is a stretched algebra of socle degree $2$ and the statement follows from Corollary \ref{Vstretched}.
Assume $e=c+3$. If $R/I$ is a stretched algebra of socle degree $3$ then we are done again by  Corollary  \ref{Vstretched}.
We may  assume that $R/I$ is not stretched. Then, its Hilbert function is
$$ HF_{R/I}\;:\qquad 1 \quad c \quad 2 \quad 0_{\longrightarrow}$$
and Proposition \ref{ineq} concludes the proof for all value of $c$, except $c=3$ or $c=4$.
However, both these cases follow from Proposition \ref{short2}.
 %Assume $c=4$.  Notice that by Remark \ref{homog}, we may assume that $R$ is a polynomial ring in $c$
% variables over a field and $I$ is a homogeneous ideal.
%First if $c=3$, write $I=(f_1, \dots, f_4,\m^3)$, where $f_i$'s are forms of degree $2$. We have to show that $\lambda(R/I^2)\,>\,24$. Notice
%$I^2=(Q ,\,f_1\m^3, \dots, f_4\m^3 ,\m^6)$ where $Q=(f_if_j\,|\,1\leq i \leq j \leq 4)$ which is  generated by $10$ quartics, say $g_1,\dots, g_{10}$. Hence
%$$I^2=(g_1,\dots, g_{10}, f_1\m^3,\dots, f_4\m^3,\m^6)\,\subseteq\,(g_1,\dots, g_{10},\m^5)$$
%and
%$$\lambda(R/I^2)\,\geq\,\lambda(R/(g_1,\dots, g_{10},\m^5))\geq 1+3+6+10+(15-10)=25\,>\,24,$$
%which implies the case $c=3$.
 %Since $4=c=\frac{1+\sqrt{1+24(t-1)}}{2}$, the proof of Proposition \ref{ineq}
%shows that we only need to avoid the case $I^2=\m^4$.
% Since $I$ is homogeneous, one could check (e.g. by direct computations using
 %CoCoA \cite{Co}) that  $I^2=\m^4$ happens only if $I=\m^2$ which gives a contradiction,
 % since $I$ has only $8$ quadratic minimal generators.
 \QED

\bigskip

%\begin{Lemma}\label{socle}
%Let $S=k[X_1,\dots,X_c]$ be a polynomial ring over a field and $q_1,\dots,q_r$ be homogeneous quadrics with $r\,\geq\,c+1$. If ${\rm ht}(\underline q)=c$ then the socle degree of $R/I$ is $\leq c-1$.
%\end{Lemma}

Next, we deal with the case $e(R/I)=c+4$. We will need the following lemma, which  shows how the Hilbert function and the associated graded ring are affected by factoring out a socle element.

\begin{Lemma}\label{socle}
Let $(R,\m)$ be a Noetherian local ring and $0\neq f\in{\rm soc}(R)$. Then $f^*$ is a non zero socle element of $R^*={\rm gr}_{\m}(R)$ and $R^*/f^*R^*\simeq {\rm gr}_{\m}(R/fR)$.
Furthermore, if $f\in \m^j\backslash \m^{j+1}$, then
${\rm HF}_{R/fR}(i)={\rm HF}_{R}(i)$ for every $i\,\neq\,j$, and ${\rm HF}_{R/fR}(j)={\rm HF}_{R}(j)-1$.
\end{Lemma}

\demo
Since $f^*\m^*\subseteq (f\m)^*=(0)$,  we deduce that $0\neq f^*$ is a socle element of ${\rm gr}_{\m}(R)$.
To prove that $R^*/f^*R^*\simeq {\rm gr}_{\m}(R/fR)$, it is enough to show the equality $(fR)^*=f^*R^*$. Clearly $f^* \in (fR)^*$ showing one inclusion.
For the other inclusion, take $b\in R$, we show that $(bf)^*\in f^*R^*$.
Indeed, if $b\notin \m$ the statement is trivial, so we may assume that $b\in \m$. But then $bf=0$ concluding the proof.
Finally, %For the part regarding the Hilbert Function
notice that for every $i$, ${\rm HF}_{R/fR}(i)={\rm HF}_{R^*/(fR)^*}(i)={\rm HF}_{R^*/f^*R^*}(i)$ and ${\rm HF}_{R}(i)={\rm HF}_{R^*}(i)$. Therefore, it is enough to prove the statement in the homogeneous settings, which is a well-known result.
\QED

\bigskip

We now show that Question \ref{Vasc2} holds true for ideals defining algebras with multiplicity $c+4$.
%For a homogeneous ring $R$ and an integer $d$, we denote by $[R]_d$ the homogeneous component of degree $d$ of $R$.

\begin{Theorem}\label{c+4}
Let $(R,\m)$ be a regular local ring containing a field of characteristic $\neq 2$ and $I$ a Cohen-Macaulay $R$-ideal that is generically a complete intersection. Assume  $e(R/I)=c+4$, where $c={\rm ecodim}(R/I)$. \\
$($a$)$ If $c\geq 5$, then $I/I^2$ $($equivalently, $R/I^2$$)$ is not Cohen-Macaulay.\\
$($b$)$ If $c\leq 4$ and $I/I^2$ $($or $R/I^2$$)$ is Cohen-Macaulay, then $R/I$ is Gorenstein.
\end{Theorem}

\demo
First of all, we employ Lemma \ref{LC} and Proposition \ref{NonDeg2} to assume that   $c\geq 3$, $I\subseteq \m^2$ and ${\rm ht}\,I=c$. Then by Proposition \ref{Generalized},\, $e(R/I^2)=(c+1)e(R/I)$.
As before, after factoring out a minimal reduction of the maximal ideal of $R/I$ generated by a system of parameters, we may assume that $I$ is an $\m$-primary ideal (not necessarily generically a complete intersection) with multiplicity (or, equivalently, length) $c+4$. %Our goal is to show that
%$\lambda(R/I^2)\,>\,(c+1)(c+4)$ under the condition that either  $c\geq 5$ or  $c\leq 4$ and $R/I$ is not Gorenstein.
Now, (a) just follows from Proposition \ref{ineq}. For (b), since $c=3$ or $ 4$, the only possible Hilbert functions of $R/I$ are:
$$\begin{array}{lllllllllllll}
(I)   & & 1 & & c & & 3 & & 0_{\longrightarrow} & & & &\\
(II)  & & 1 & & c & & 2 & & 1 & & 0_{\longrightarrow} & &\\
(III) & & 1 & & c & & 1 & & 1 & & 1 & & 0_{\longrightarrow}
\end{array}$$

Since in case (I) the algebra is short, Proposition \ref{short2} completes the proof. In case (III), the algebra is stretched, hence the conclusion follows by Theorem \ref{Tstretched}.
We are then left with case (II). Assume by contradiction that $R/I$ is not Gorenstein. We show that this forces $\lambda(R/I^2)\,>\,(c+1)(c+4)$, contradicting the previously proved equality $e(R/I^2)=(c+1)e(R/I)$.

Since $R/I$ is not Gorenstein,  $R/I$  has a socle element of initial degree $1$ or $2$. We first deal with the case where there exists a  socle element of  initial degree  $1$. Since it is part of a minimal generating set of the maximal ideal, we may write $\m=(x_1, \ldots, x_c)$, where $x_1$ is a socle element of $R/I$. Notice that, after passing to the completion, we may assume that $R$ is a power series ring over $k$ in the variables $x_1,\dots,x_c$.

Let $c=4$. Then  $I=(x_1\m,q_1,\dots,q_4, g,\m^4)$, where $g$ is a combination (with unit coefficients) of monomials of degree $3$ in $x_2,x_3,x_4$ and all the $q_i's$ have initial degree 2 and involve only $x_2, x_3, x_4$.
Since $x_1$ is a socle element, by Lemma \ref{socle} the Hilbert function of $R/(I, x_1)$ is $$1 \qquad 3 \qquad 2 \qquad 1 \qquad 0_{\longrightarrow}\,{}$$

Hence
   %if one sets $\n=(x_2,x_3,x_4)$ then $(q_1,\dots,q_4)\n\neq \n^3$. Hence
   $(x_1\m^2, (q_1,\dots,q_4)\m)\neq\m^3$ and
$(x_1^2\m^2, x_1(q_1,\dots,q_4)\m)$ is strictly contained in $x_1\m^3$.
 Now $I^2\subseteq (x_1^2\m^2, x_1(q_1,\dots,q_4)\m, q_iq_j, 1\leq i \leq j \leq 4, \m^5)\subsetneq (x_1\m^3,  q_iq_j, 1\leq i \leq j \leq 4, \m^5)$,
\begin{eqnarray*}
\lambda(R/I^2)&> &\lambda (R/ (x_1\m^3,  q_iq_j, 1\leq i \leq j \leq 4, \m^5))\\
&\geq & 1 + 4 + 10 + 20 + 5 = 40=(c+1)e(R/I),
\end{eqnarray*}
and we are done.

Now, assume $c=3$. Since $(c+1)e(R/I)=28$, our goal is to show that $\lambda(R/I^2)\geq 29$.
Similarly to the above argument,  $I=(x_1\m, q, g,\m^4)$, where $q$  and $g$ involve only $x_2, x_3$, $q$ has initial degree 2 and $g$  is a combination (with unit coefficients) of monomials of degree $3$ in $x_2,x_3$.
Then $I^2\subseteq(x_1^2\m^2,x_1 q\m, q^2, x_1g\m, qg, \m^6)$ and similar to the above one gets that $$\lambda(R/I^2)\geq \lambda(R/(x_1^2\m^2,x_1 q\m, q^2, x_1g\m, qg, \m^6))\geq 1+3+6+10+6+3=29.$$
%$HF_{R/I^2}(4)\geq 6$, so to finish the proof, it is enough to show that $HF_{R/I^2}(5)\geq 3$.
%Notice that $(I,x_1)=(x_1, q, g,\m^4)$, hence $(I^2, x_1)=(x_1, q^2, qg, g^2, q\m^4, g\m^4,\m^8)$.  Clearly we have
%$HF_{R/(I^2,x_1)}(5)\geq 3$, proving that $HF_{R/I^2}(5)\geq 3$ and then finishing the proof in this case.
Hence, from now on we may assume that $R/I$ has no socle elements in $\m\backslash \m^2$.

First, consider $c=4$. Then, there exists an element $q_9\in \m^2 \backslash \m^3$ with $q_9\in {\rm soc}(R/I)$ and $q_9^*\notin I^*$. By Lemma \ref{socle}, $K:=(I,q_9)$ defines a stretched Artinian  algebra of socle degree 3 and there exist $x_1,\dots,x_4$ in $\m\backslash \m^2$ so that a generating set of $K$ can be described as in Theorem \ref{SS}. Observe that $I^2\subseteq L=(x_1, x_2, x_3)H+(x_1, x_2, x_3)x_4^4+ (x_4^6)$, where  $H$ is as in Proposition \ref{Isquare}. Since
$$
e(R/I^2)\!=\!\lambda(R/I^2)\geq \lambda (R/L)\geq 1 + 4 + \binom{4+1}{2} + \binom{4+2}{3} + 4 + 1=40=(c+1)e(R/I),
$$
to finish the proof, it is enough to prove that $I^2\subsetneq L$. If $R/K$ has type $3$ or $4$, we are done by Proposition \ref{Isquare}. So assume $R/K$ has type $1$ or $2$, then $K$ has nine minimal generators, say, $q_1, \ldots, q_9$. The same idea of Remark \ref{homog} shows that we can further assume that $R$ is $k[x_1,\dots,x_4]$ and $I$ is generated by polynomials whose homogeneous components have degrees $2$ and $3$ (this follows from $\m^4\subseteq I\subseteq \m^2$ and the fact that $k[x_1,\dots,x_4]/I \simeq k[\![x_1,\dots,x_4]\!]/I$).

Notice that $I=(q_1, \ldots, q_8, q_9\m)$. Also, the description of $K$ yields that for every $i$ one has $q_i=Q_i + a_ix_4^3$ for some $Q_i$ homogeneous of degree $2$ and $a_i$ (either zero or unit) in $R$.
Since $I^2=(q_iq_j,  q_iq_9\m, 1\leq i, j\leq 8, q_9^2\m^2)\subseteq K^2=(q_iq_j,  1\leq i, j\leq 9)\subseteq L$, if $I^2=K^2=L$, an application of Nakayama's Lemma -- after going modulo $(q_iq_j,   1\leq i, j\leq 8)$ -- shows that  $L=(q_iq_j,   1\leq i, j\leq 8)=(q_1, \ldots, q_8)^2$.% (this is an application of Nakayama's Lemma after going modulo $(q_iq_j, 1\leq i \leq j \leq 8)$) . % If we denote by $^{-\!\!-}$ the image modulo $(q_iq_j, 1\leq i \leq j \leq 8)$ we have that $\overline{K}^2=(q_iq_9, 1\leq i \leq 9)=\overline{I}^2=(q_iq_9\m, 1\leq i \leq 8, q_9^2\m^2)\subseteq \overline{\m}\overline{K}^2$, thus showing that $\overline{K}^2=\overline{(0)}$

Set $\tilde{R}=R[t]$, where $t$ is a variable over $R$. For $1\leq i\leq 9$,
let $\tilde{q}_i$ be  obtained from $q_i$ by  replacing $x_4^3$ by $x_4^2 t$, and set $\tilde{K}=(\tilde{q}_1, \ldots, \tilde{q}_9)\tilde{R}$ and $\tilde{J}=(\tilde{q}_1, \ldots, \tilde{q}_8)\tilde{R}$. It is easy to see that
$\tilde{K}^2=\tilde{L}=(x_1, x_2, x_3)H+(x_1, x_2, x_3)x_4^3t+ (x_4^4 t^2)$. Now we have $\tilde{J}^2\subseteq \tilde{L}$.

Let $\pi(\cdot)$ denote the images modulo $t-x_4$, then $\pi(\tilde{J}^2)=\pi(\tilde{L})=L$. Since $\tilde{L}$ is a monomial ideal, one can check that $t-x_4$ is regular on $\tilde{R}/\tilde{L}$. Hence
$\tilde{L}=\tilde{J}^2+(t-x_4)\tilde{R}\cap \tilde{L}=\tilde{J}^2+(t-x_4)\tilde{L}$. By Nakayama's lemma, $\tilde{L}=\tilde{J}^2$. Now set $t=1$. Then, $\tilde{q}_1, \ldots, \tilde{q}_8$ become homogeneous of degree two, say $q_1^{\prime}, \ldots q_8^{\prime}$ and we have that $(q_1^{\prime}, \ldots q_8^{\prime})^2=\m^2$. This fact gives a contradiction to \cite[Theorem~2.4]{COR} (see also the proof of Proposition \ref{short2}).

Assume, then, that $c=3$.
Observe that $I=(g_1, g_2, g_3, g_4, g_5\m)$,  $(I,g_5)=(g_1, g_2, g_3, g_4, g_5)$, where all the $g_i's$ have  initial degree 2 and $g_5\in {\rm Soc}(R/I)$, $g_5 \notin I$. There are two cases. First case, $(I, g_5)$ is a stretched Gorenstein ideal. Then, $I^2=(g_ig_j, 1\leq i,j\leq 4, g_ig_5\m, 1\leq i\leq 4, g_5^2\m^2)\subseteq (I, g_5)^2=(g_ig_j, 1\leq i,j\leq 5)= (x_1,x_2)H+(x_1, x_2)x_3^4+(x_3^6)$. Since the number of minimal generators of $(I, g_5)^2$ is $15$, it is generated minimally by $g_ig_j, 1\leq i, j\leq 5$. The ideal $I^2$ contains only $10$ minimal generators of them, hence $\lambda((I,g_5)^2/I^2)\geq 5$. Therefore $\lambda(R/I^2)= \lambda(R/(I, g_5)^2)+\lambda((I,g_5)^2/I^2)\geq (1 + 3 + \binom{3+1}{2} + \binom{3+2}{3} + 3 + 1)+5=29> 28=(c+1)e(R/I)$.

Second case, $R/(I,g_5)$ is not Gorenstein. If $\tau(R/(I, g_5))=2$, then $(I,g_5)^2=(x_1, x_2)H+(x_2x_3^4-ux_2^3x_3, x_3^6+u^2x_2^4)$ which has 13 minimal generators. $I^2$ contains at most 10 of them. Hence,
$$\lambda(R/I^2)= \lambda(R/L)+\lambda(L/(I, g_5)^2)+\lambda((I,g_5)^2/I^2)\geq 24+2+3=29,$$
and the proof is finished.
Finally, if $\tau(R/(I,g_5))=3$, then $(I,g_5)=(x_1, x_2)\m+(x_3^4)$. Hence $\lambda(R/(I,g_5)^2)= 1+3+6+10+3+3+1+1=28$.
Since $(I,g_5)$ has 12 minimal generators, then $\lambda(R/I^2)\geq 28+2=30$, providing the desired contradiction.
\QED

\medskip

As  a consequence of  Theorems \ref{c+3} and \ref{c+4}, we give a positive answer to Question \ref{Vasc2} for ideals defining algebras with low multiplicity.

\begin{Corollary}\label{LowMult}
 Question \ref{Vasc2} holds true for $I$ with  $e(R/I)\,\leq\,c+4$, where $c={\rm ecodim}(R/I)$.
\end{Corollary}

In particular, Question \ref{Vasc2} holds true for some irreducible algebroid curves.
\begin{Corollary}
Let $S=k[\![f_1(t),\ldots,f_n(t)]\!]$ for some $f_1,\dots,f_n$ in $k[\![t]\!]$ and assume none of the $f_i$'s is redundant.
Let $R=k[\![X_1,\ldots,X_n]\!]$ and $\p$ the prime $R$-ideal defining $S$. If $a_i\!:={\rm indeg}(f_i)\leq n+3$ for some $i$, then Question \ref{Vasc2} holds true for $\p$.\\
In particular, if $S=k[\![t^{a_1},\ldots,t^{a_n}]\!]\simeq R/\p$, where none of the $t^{a_i}$ is redundant, and $a_i \leq n+3$ for some $i$, then Question \ref{Vasc2} holds true for $\p$.
\end{Corollary}

\demo
In this setting, it is well-known that $n=c+1$ and $e(R/\p)\leq a_j$ for every $j$. Hence  $a_i\leq n+3$ for some $i$ implies that $e(R/\p)\leq c+4$. Now apply Corollary \ref{LowMult}.
\QED

%\medskip

\section{Counterexamples and sharpness of the main result}

In this section we provide counterexamples to Question \ref{Vasc2} showing that, in general, it has a negative answer. Also, these examples shows the sharpness of Theorem \ref{Main} (d)-(e).%We show at the same time that the results proved in the previous sections (in particular \ref{LowMult}) are actually sharp.

We first exhibit an example found by  J. C. Migliore using CoCoA \cite{Co}.
Then, we  deform it using the theory of universal linkage (see \ref{prime}) to obtain a prime ideal that is also a counterexample to Question \ref{Vasc2}. Finally, we conjecture that in $\mathbb P^c$  a certain number of points (function of $c$) gives a counterexample to Question \ref{Vasc2} in {\em any} embedding codimension $c\geq 5$.

\begin{Example}\label{count}
$($a$)$ Let $S=k[a,b,c,d,e,f]$ be a polynomial ring over $k$, where $k$ is either $\mathbb Q$ or $\mathbb Z/31991\mathbb Z$ and let $M=(a,
dots,f)$. Then a homogeneous level ideal $I$ in $S$ defining $10$ general points gives a negative answer to the homogeneous version of Question \ref{Vasc2}, i.e., $I/I^2$ and $S/I^2$ are Cohen-Macaulay but $S/I$ is not Gorenstein. %$($in fact its Cohen-Macaulay type is $4$$)$. %that is generically a complete intersection, with both $I$ and $I/I^2$ Cohen-Macaulay but $R/I$ is not Gorenstein.
Also, $S/I$ is a short algebra with ${\rm socdeg}(S/I)=2$, $\tau(S/I)=4$, $c=5$ and $e(S/I)=c+5$.

$($b$)$ There exist a regular local ring $(R,\m)$ and a prime ideal $\p$ such that $\p/\p^2$ and $R/\p^2$ are Cohen-Macaulay but $R/\p$ is not Gorenstein. $R/\p$ has multiplicity $c+5$, where $c={\rm ecodim}(R/\p)$.
\end{Example}

\demo
(a) In  $S=k[a,b,c,d,e,f]$,
 using the package `IdealOfProjectivePoints' (see \cite{ABKR}), one obtains the following homogeneous ideal which defines a set of $10$ general points in $\mathbb P^5$:
 $$
  \begin{array}{l}
 I=( e^2f + 2963bf^2 + 4964cf^2 + 5333df^2 - 13261ef^2, \\
     af - 13894bf + 12842cf + 4036df - 2985ef,
     de + 3056bf + 12160cf + 971df + 15803ef, \\
     ce - 2357bf - 14460cf + 3040df + 13776ef,
     be - 8504bf + 1159cf - 1581df + 8925ef, \\
     ae - 9147bf + 1379cf + 4167df + 3600ef,
     cd + 7380bf + 5885cf + 6255df + 12470ef, \\
     bd - 11676bf - 2833cf - 13277df - 4206ef,
     ad + 5555bf + 2017cf + 2100df - 9673ef, \\
     bc + 5653bf - 6596cf - 8208df + 9150ef,
     ac - 2335bf - 10387cf + 514df + 12207ef, \\
     ab - 8324bf - 7688cf - 4252df - 11728ef,
     b^2f + 15536bf^2 + 1265cf^2 + 9888df^2 + 5301ef^2, \\
     d^2f + 10625bf^2 - 11725cf^2 + 9514df^2 - 8415ef^2, \\
     c^2f + 11390bf^2 + 7112cf^2 - 10319df^2 - 8184ef^2).
  \end{array}$$
%Now, since $I$ defines a set of projective points then ${\rm dim}(S/I)=1$.
%Furthermore -- u
Using CoCoA for instance, one can see that the H-vector of $S/I$ is
$$1\qquad 5 \qquad 4$$
 proving that $c=5$, $e(S/I)=c+5$, and $S/I$ is a short algebra with   ${\rm socdeg}(S/I)=2$.
By the  Betti diagram of $S/I$, one has that $S/I$ is level which implies that $\tau(S/I)=4.$
Finally the Betti diagram of $S/I^2$ shows that $S/I^2$ is Cohen-Macaulay, yielding that $I/I^2$ is Cohen-Macaulay too.
%Use R::=QQ[a,b,c,d,e,f];
%I:=IdealOfProjectivePoints(GenericPoints(10));
%HVector(R/I);
%
%[1, 5, 4]
%-------------------------------
%J:=I^2;
%HVector(R/J);
%
%[1, 5, 15, 35, 4]
%-------------------------------
%BettiDiagram(R/J);

%        0    1    2    3    4    5
%-----------------------------------
% 0:     1    -    -    -    -    -
% 1:     -    -    -    -    -    -
% 2:     -    -    -    -    -    -
% 3:     -   66  204  240  120   15
% 4:     -    -    -    -    -    4
%-----------------------------------
%Tot:    1   66  204  240  120   19
%-------------------------------
%BettiDiagram(R/I);
%        0    1    2    3    4    5
%-----------------------------------
% 0:     1    -    -    -    -    -
% 1:     -   11   20    5    -    -
% 2:     -    -    -   16   15    4
%-----------------------------------
%Tot:    1   11   20   21   15    4
%-------------------------------

(b) The first statement follows by part (a) and Proposition \ref{prime}. The second statement follows from the fact that $e(R/\p)\leq e(S/J)=c+5$, because $(R,\p)$ is a deformation of a $(S_M,J_M)$. However, if $e(R/\p) < c+5 $ it would contradict Theorem \ref{LowMult}. Therefore, $e(R/\p)=c+5$.
%Let $I$ be the homogeneous ideal of part (a) and write $I_{\m}$ for its localization at the homogeneous maximal ideal. Clearly $I_{\m}$ is a counterexample to Vasconcelos Conjecture, in particular it is still generically a complete intersection, is Cohen-Macaulay and its conormal module is Cohen-Macaulay.\\
%Let $L_2(I)\subseteq S[Z_1,\dots,Z_m]$ be the second generic link of $I$ in a polynomial extension of $I$. By Theory of Linkage (see for instance \cite{HU2}) there exist elements $a_1,\dots,a_m\in S$ with $Z_1-a_1,\dots,Z_m-a_m$ a regular sequence on $L_2(I)R$ where $R=S[Z_1,\dots,Z_m]_{(\m,\underline Z - \underline a)}$ and $(L_2(I),\underline Z - \underline a)/(\underline Z - \underline a)\simeq I$ in $R/(\underline Z - \underline a)\simeq S$.\\ %(in other words, there exists a localization $R$ of a polynomial ring over $S$ such that $(U,L_2(I))$ is a deformation of $(S,I)$).\\
%Since $(R,L_2(I)R)$ is a deformation of $(S,I)$ we have that both $L_2(I)R$ and $L_2(I)R/(L_2(I)R)^2$ are Cohen-Macaulay too, $\tau(R/L_2(I))=4$. Since $R$ is a regular local ring every prime ideal is generically a complete intersection, hence to finish is enough to show that $L_2(I)R$ is a prime ideal of $R$. However since $I$ is generically a complete intersection by Theory of Linkage (see for instance \cite{HU}) $L_2(I)$ is a prime ideal in $S[Z_1,\dots,Z_m]$ and hence $L_2(I)R$ is prime ideal in $R$.
\QED
\bigskip

%Juan C. Migliore, using CoCoA produced an example of a set of $10$ generic points in $\mathbb P^5$ such that both $R/I$ and $I/I^2$ are Cohen-Macaulay, but $R/I$ is not Gorenstein, as can be seen by the Hilbert Function $1 \qquad 5 \qquad 4 \qquad 0_{\longrightarrow}$. One should also remark that the minimal graded free resolution for such example is pure.\\
%This is the first counterexample of Vasconcelos Conjecture that we are aware of. Furthermore, the next section of the paper shows that - in some sense - this counterexample is the smallest possible. In fact, we show that if an algebra has multiplicity at most $c+4$ then Vasconcelos Conjecture holds.
%Hence a counterexample should have multiplicity at least $c+5$. And $10$ points in $\mathbb P^5$ actually have multiplicity $c+5$.\\
Since both these examples have multiplicity $c+5$, they shows that  Corollary \ref{LowMult} is sharp.
Furthermore, Example \ref{count} shows that the numerical estimates presented in Proposition \ref{short2} are very sharp. Indeed, for the case $c=5$,  Proposition \ref{short2} implies that  $I/I^2$ is not Cohen-Macaualy for any $q\,\neq\,11$ (using the notation of Proposition \ref{short2}). While for $q=11$, we have the above counterexample, showing that $I/I^2$ is Cohen-Macaulay.

Similarly, when $c=6$, Proposition \ref{short2} (a) and (b) show that $I/I^2$ is not Cohen-Macaulay for %$q>16$ and $q\leq 15$. Hence, $I/I^2$ is not Cohen-Macaulay when 
any $q\neq 16$ . However, for $q=16$, J. C. Migliore found with CoCoA \cite{Co} a set of $12$ general points in $\mathbb P^6$ with $I/I^2$ Cohen-Macaulay but $R/I$ not Gorenstein.

It would be natural to  conjecture that a set consisting of $2c$ general points in $\mathbb P^c$ gives a negative answer to Question \ref{Vasc2} for any $c\geq 5$. However, this is false. Indeed, a set of $14$ general points in $\mathbb P^7$ does {\em not} give a counterexample to the question.

Instead, we conjecture that the estimate of Proposition \ref{short2} is sharp and if we go past this number we get counterexamples to Question \ref{Vasc2} in any codimension at least $5$. We now state this more precisely. Let $\left\lceil x \right\rceil$ denote the smallest integer bigger than or equal to a given real number $x$. Then,
\begin{Conjecture}\label{counterex}
For any $c\geq 5$, the homogeneous ideal $I$ defined by a set of $1+c+\lceil\frac{c(c-1)}{6}\rceil$ general points in $\mathbb P^c$ has the property that $I/I^2$ are Cohen-Macaulay, but $R/I$ is not Gorenstein.
\end{Conjecture}

%Notice that if Question \ref{counterex} had a positive answer, we could use it together with Lemma \ref{prime} to find for any $c$ even a prime ideal in a $c$-dimensional regular local ring that is a counterexample to Question \ref{Vasc2}.
Our conjecture is supported by data obtained by J. C. Migliore. In fact, using CoCoA, he checked that, for any $c\leq 9$, Conjecture \ref{counterex} holds true and indeed provides the smallest set of (general) points in $\mathbb P^c$ giving a counterexample to Question \ref{Vasc2}. %However, for $c\geq 10$, we have no data, since the computations needed to check Question \ref{counterex} become too expensive. As of now, there are no known counterexamples to Question \ref{counterex}.

\bigskip

{\bf Acknowledgments.}
The authors warmly thank the AMS, for both organizing and generously supporting the MRC week in the wonderful
location of Snowbird, UT; the inspiration for many ideas of this paper came out during that week.
We are very grateful to Claudia Polini -- for giving us precious suggestions and ideas, and contaging us with her enthusiasm, and Bernd Ulrich -- for suggesting several improvements and a careful reading of the present paper.
We are also in debt with Craig Huneke, for several insightful remarks given during the MRC week, and J. C. Migliore, whose experience and help produced the homogeneous counterexamples to Question \ref{Vasc2}.

Least, but no last, we would like to thank warmly the organizers of the 2010 MRC week, D. Eisenbud, C. Huneke, M. Musta\c t\u a and C. Polini.
Their strong support and passionate work made the MRC week an unforgettable experience.
%
%-------------------
%
%OLD
%We are indebted to Claudia Polini for generously giving us precious suggestions and ideas
%and contaging us with her enthusiasm. We are very grateful to Bernd Ulrich for suggesting several improvements and for a careful reading of the present paper.
%We are also in debt with Craig Huneke, for several insightful remarks he gave us during the MRC week, and J. C. Migliore, whose experience and computations with CoCoA have been extremely helpful.
%Least but no last, we would like to thank warmly the organizers of the MRC, D. Eisenbud, C. Huneke, M. Musta\c t\u a and C. Polini.
%Their strong support and passionate work made the MRC week a successful and unforgettable experience.

\smallskip

%\end{document}

%{\small

\end{document}